\documentclass[11pt,a4paper]{article}
\usepackage[a4paper,top=2.56cm,bottom=2.56cm,left=2.56cm,right=2.56cm,marginparwidth=1.75cm]{geometry}
\usepackage{amsmath,amsfonts,amssymb,amsthm}
\usepackage{mathtools}
\usepackage{graphicx}
\usepackage[colorlinks=true, allcolors=blue]{hyperref}

\numberwithin{equation}{section}
\newtheorem{theorem}{Theorem}[section]
\newtheorem{lemma}[theorem]{Lemma}
\newtheorem{proposition}[theorem]{Proposition}
\newtheorem{corollary}[theorem]{Corollary}
\theoremstyle{definition}
\newtheorem{definition}[theorem]{Definition}
\newtheorem{example}[theorem]{Example}
\theoremstyle{remark}
\newtheorem{remark}[theorem]{Remark}
\theoremstyle{plain}
\makeatletter
\def\@fnsymbol#1{\ifcase#1\or a\or b\or c\or d\fi}
\newcommand{\TITLE}[1]{\def\@TITLE{#1}}
\newcommand{\SHORTTITLE}[1]{}
\newcommand{\AUTHORS}[1]{\def\@AUTHORS{#1}}
\newcommand{\KEYWORDS}[1]{\def\@KEYWORDS{#1}}
\newcommand{\AMSSUBJ}[1]{\def\@AMSSUBJ{#1}}
\newcommand{\AMSSUBJSECONDARY}[1]{\def\@AMSSUBJSECONDARY{#1}}
\newcommand{\ACCEPTED}[1]{}
\newcommand{\ABSTRACT}[1]{\def\@ABSTRACT{#1}}
\newcommand{\FIRSTPAGE}{%
  \title{\@TITLE}%
  \author{\@AUTHORS}%
  \date{}%
  \maketitle
  \begin{abstract}%
    \noindent\@ABSTRACT\par\vskip 1em\relax
    {\footnotesize
      \noindent\textbf{Keywords: }\@KEYWORDS.\par
      \noindent\textbf{AMS MSC 2020: }Primary \@AMSSUBJ, Secondary \@AMSSUBJSECONDARY.\par
    }%
  \end{abstract}%
  \hypersetup{%
    pdftitle={Hitting times of non-reversible Markov chains via the singular value decomposition},
    pdfauthor={Thomas M. van Belle, Anton Klimovsky},
    pdfsubject={Probability Theory}}%
}
\AtBeginDocument{\FIRSTPAGE}
\makeatother

\usepackage[english]{babel}
\usepackage[utf8]{inputenc}
\usepackage[T1]{fontenc}
\usepackage[svgnames]{xcolor}
\usepackage{bbm}
\usepackage{mathrsfs}
\usepackage{caption}
\usepackage{subcaption}
\usepackage{float}
\usepackage{url}
\usepackage{enumitem}
\usepackage{pgfplots}
\pgfplotsset{compat=1.16} %
\usepackage{pgfplotstable}

\definecolor{light-gray}{gray}{0.97}
\definecolor{Blue}{rgb}{0,0,0.5}
\definecolor{Green}{rgb}{0,0.75,0.0}
\definecolor{LightGray}{rgb}{0.6,0.6,0.6}
\definecolor{DarkGray}{rgb}{0.3,0.3,0.3}

\newcommand{\defeq}{\mathrel{\mathop:}=}
\newcommand{\eqdef}{=\mathrel{\mathop:}}

\def\<{\langle}
\def\>{\rangle}

\def\R{{\mathbb{R}}}  
\def\N{{\mathbb{N}}}  
\def\P{{\mathbb{P}}}  
\def\Z{{\mathbb{Z}}}  
\def\C{{\mathbb{C}}} 
\def\E{{\mathbb{E}}}

\let\cal=\mathcal

\def \cF {{\cal F}}
\def \cH {{\cal H}}

\def \cU {{\cal U}}

\usepackage[
backend=biber,
style=alphabetic,
giveninits=true,
sortlocale=auto,
date=year,
natbib=true,
url=false,
doi=true,
isbn=false,
eprint=false,
maxnames=999,
maxbibnames=999,
maxcitenames=50
]{biblatex}
\appto{\bibsetup}{\sloppy}
\renewbibmacro{in:}{\ifentrytype{article}{}{\printtext{\bibstring{in}\intitlepunct}}}

\TITLE{Hitting times of non-reversible Markov chains via the singular value decomposition}
\SHORTTITLE{Hitting times of non-reversible Markov chains via the SVD}
\AUTHORS{Thomas\,M. van Belle\thanks{Work carried out at the Fakult\"at f\"ur
Mathematik, Universit\"at Duisburg-Essen, Thea-Leymann-Str.~9, 45127 Essen,
Germany.} \and Anton
Klimovsky\thanks{Julius-Maximilians-Universit\"at W\"urzburg, Institut f\"ur
Mathematik, Lehrstuhl f\"ur Mathematik VIII: Angewandte Stochastik,
Emil-Fischer-Str.~30, 97074 W\"urzburg, Germany.
\texttt{anton.klimovsky@mathematik.uni-wuerzburg.de}}}

\AMSSUBJ{60J10}
\AMSSUBJSECONDARY{60B20, 05C80}
\ACCEPTED{--}
\ABSTRACT{We develop a singular value approach to hitting times for finite non-reversible
Markov chains. For an irreducible chain with transition matrix $P$, generator
$L=I-P$ and stationary distribution $\pi$, we write the generalised inverse
formula for the hitting times $H_{ij}$ in the singular values and singular vectors
of $L$. This gives an exact bilinear expansion in the left and right singular vector
bases. It needs no
reversibility, normality or diagonalisability, it reduces to the
Lov\'asz spectral formula in the reversible case, and its truncation to the slow
singular modes has an error that involves no eigenvector condition number.

For simple random walk on the directed Erd\H{o}s--R\'enyi graph with
$\log(n)\ll d=np_n$, we prove laws of large numbers with explicit exponential
error bounds: Kemeny's constant satisfies $H^i=n(1+o(1))$, and
$\pi_jH_{ij}=1+o(1)$ uniformly over ordered pairs $i\neq j$. If $d\geq n^{\delta}$
for some $\delta>0$, then $H^i=n-2+o(1)$ in probability. If $\limsup_n p_n<1$, the
empirical distribution of the centred singular values
$\sqrt{2d/(1-p_n)}\,(\sigma_k(L)-1)$ converges weakly in probability to the Wigner
semicircle law. Although the generator is non-normal and has a complex spectrum,
the leading fluctuations of $L^tL-I$ are governed by the Hermitian part of the
centred transition matrix, giving rise to the semicircle law.
}

\begin{document}

\section{Introduction}\label{sec:introduction}
Let $(X_t)_{t \geq 0}$ be an irreducible Markov chain on a finite state space $[n]$
with transition matrix $P$ and stationary distribution $\pi$, and let
\begin{equation}
	H_{ij} \defeq \E_i\left[\inf\{t \geq 0 : X_t = j\}\right], \qquad L \defeq I-P.
\end{equation}
For a \emph{reversible} chain, $H_{ij}$ has the spectral representation of
Lov\'asz~\cite{lovasz_random_1996} (Theorem~\ref{thm:lovaszhittingtimes}): a sum
over the eigenmodes of $P$, weighted by the inverse spectral gaps
$(1-\lambda_k)^{-1}$. It underlies much of the literature on hitting times of
random walks on random
graphs~\cite{lowe_hitting_2014, helali2019hitting,lowe2024spectral,lowe2025central,ottolini2024concentration}.

For the simple random walk on a directed graph, reversibility fails. Kemeny's
constant still equals $\sum_{k\geq2}(1-\lambda_k)^{-1}$ for every irreducible
chain, with the eigenvalues counted with algebraic multiplicity
(Remark~\ref{rem:beyondreversible}), but the expansion of the individual
$H_{ij}$ is lost. Left and right eigenvectors differ and are not orthogonal, so an
eigenvalue expansion may be badly conditioned. Worse, $P$ need not be
diagonalisable, and we give an irreducible three-state chain without an
eigenbasis (Example~\ref{ex:defective}).

The singular value decomposition of $L$ always exists. Its singular values are real and
nonnegative, and the left and right singular vectors form two orthonormal bases.
Chatterjee~\cite{chatterjee_spectral_2023} used the singular values of the
similarity-transformed generator
\begin{equation}\label{eq:Lhatintro}
	\hat L \defeq \operatorname{diag}(\pi)^{1/2}\,L\,\operatorname{diag}(\pi)^{-1/2}
\end{equation}
to define a relaxation time for non-reversible chains; see
also~\cite{xu2026singularvaluegap}. In this sense the smallest non-zero singular
value of $\hat L$ is a quantitative mixing condition. Singular values are,
however, not eigenvalues of $P$. Hitting times have a classical expression through
a generalised inverse of $L$~\cite{hunter_generalized_1982,meyer_group_1975}. We
choose the Moore--Penrose inverse and write this expression in the singular
system of $L$. The result is an exact expansion in two orthonormal bases. It holds
for every irreducible chain, defective ones included, and its truncation to the
slow modes has an error that involves no eigenvector condition number. Meeting
times are treated by the same route in~\cite{vanbelle2024meeting}.

Our contributions are the following.
\begin{enumerate}[label=(\arabic*)]
	\item An exact representation of the classical generalised inverse formula for
	$H_{ij}$ in the singular values and singular vectors of $L$, for every
	irreducible chain on a finite state space (Theorem~\ref{thm:Hittingsingular}). Its truncation to the modes with the
	smallest non-zero singular values has an error bound with no eigenvector
	condition number (Proposition~\ref{prop:slowmodes}).
	\item Laws of large numbers for the hitting times of the simple random walk on
	the directed Erd\H{o}s--R\'enyi graph $\overrightarrow{G}(n,p_n)$ with
	$\log(n)\ll d = np_n$ (Theorem~\ref{thm:hittingdirectedER}), including the
	second-order law $H^i = n-2+o(1)$ for Kemeny's constant.
	\item A semicircle law for the singular values of the generator of this walk
	(Theorem~\ref{thm:semicircle}).
\end{enumerate}

\subsection{Results}

\begin{theorem}[Singular decomposition of hitting times]
	\label{thm:Hittingsingular}
	Let $P$ be the transition matrix of an irreducible Markov chain on the finite state space $[n]$, with stationary distribution $\pi$. Denote by $(\sigma_{k})_{k\in[n]}$ the singular values of $L = I - P$, and denote the left and right singular vectors by $(u_k)_{k\in[n]}$ and $(v_k)_{k\in[n]}$ respectively. Then, the expected hitting time from state $i$ to state $j$ is given by
	\begin{equation}\label{eq:Hittingsingular}
		H_{ij} = \sum_{k=1}^{n-1} \frac{1}{\sigma_k} \left(v_{ki} - v_{kj}\right)\left(u_k^t\mathbf{1}_n- \frac{u_{kj}}{\pi_j}\right),
	\end{equation}
	where $\mathbf{1}_n$ denotes the vector of ones of length $n$.
\end{theorem}

For an irreducible chain the eigenvalue $1$ of $P$ is simple, so
$\operatorname{rank}(L) = n-1$ and $\sigma_{n-1}(L)>0$. The proof
(Section~\ref{sec:hittingtimessvd}) follows the reversible argument, with the
Moore--Penrose inverse $L^{+} = \sum_{k<n}\sigma_k^{-1}v_ku_k^t$ in place of the
fundamental matrix. In terms of $L^+$ the identity reads
$H_{ij} = (L^{+}\mathbf{1}_n)_i - (L^{+}\mathbf{1}_n)_j + ((L^{+})_{jj} - (L^{+})_{ij})/\pi_j$
(equation~\eqref{eq:HijviaMP}). This is the Moore--Penrose case of Hunter's
formula for mean first passage times through an arbitrary generalised inverse of
$L$~\cite{hunter_generalized_1982}, and its group-inverse counterpart is due to
Meyer~\cite{meyer_group_1975}. Since $L^{+}$ is unique, the right-hand side
of~\eqref{eq:Hittingsingular} does not depend on the choice of the singular vectors.

In the Lov\'asz formula the expansion runs over one orthonormal eigenbasis.
Theorem~\ref{thm:Hittingsingular} replaces it by two orthonormal bases, the right
and the left singular vectors, and the expansion is bilinear in them. Without reversibility the eigenvectors of $P$ need not be orthogonal,
and they need not form a basis (Example~\ref{ex:defective}). The left and right
singular vectors of $L$ always form two orthonormal bases, and
Theorem~\ref{thm:Hittingsingular} expands $H_{ij}$ exactly in these modes.
Proposition~\ref{prop:slowmodes} truncates the expansion to the slow modes with an
error that involves no eigenvector condition number. For the product-form lumped
chains of Example~\ref{ex:lumped} the block-level and within-block contributions
separate exactly. On the directed
Erd\H{o}s--R\'enyi graph, the location of the singular values of $L$ and an
alignment bound for the left and right singular vectors give the first order of
Kemeny's constant through Theorem~\ref{thm:Hittingsingular}
(Appendix~\ref{sec:svdroute}). Theorem~\ref{thm:semicircle} identifies the limiting
empirical distribution of the centred and rescaled singular values.

Applied to $\hat L$, the formula
returns the Lov\'asz formula~\eqref{eq:Hittingeigen} in the reversible case
(Remark~\ref{rem:Ltilde}). The two versions are both exact. The singular values of
$\hat L$ lie in $[0,2]$, whereas those of $L$ carry the degree fluctuations and
can be unbounded (on a star, $\sigma_1(L)$ grows like $\sqrt n$). On the oriented
cycle the formula gives $H^i = (n-1)/2$, against $(n^2-1)/6$ on the unoriented
cycle (Remark~\ref{rem:cyclecomparison} and
Appendix~\ref{sec:circulantexample}).

Our main non-reversible example is the simple random walk on the directed
Erd\H{o}s--R\'enyi graph $\overrightarrow{G}(n,p_n)$ of
Definition~\ref{def:directedER}. We write $d = np_n$, and
$H^i \defeq \sum_{j=1}^n \pi_j H_{ij}$ and $H_j \defeq \sum_{i=1}^n \pi_i H_{ij}$
for the average target and average starting hitting times.

\begin{theorem}[Laws of large numbers for hitting times on $\protect\overrightarrow{G}(n,p_n)$]
	\label{thm:hittingdirectedER}
	Consider simple random walk on $\overrightarrow{G}(n,p_n)$ with $\log(n) \ll d$, and let $\overrightarrow{\P}_{n,p}$ denote the law of the model.
	\begin{enumerate}[label=\textup{(\alph*)}, ref=\textup{(\alph*)}]
		\item\label{it:Hitarget} For every $\epsilon \in (0,1]$ and all $n \geq 300/\epsilon$,
		\begin{equation}\label{eq:meanhittingdirectedER}
			\overrightarrow{\P}_{n,p}\left(\left|\frac{1}{n}H^i - 1\right| > \epsilon\right) \leq 4n\exp\left(-\frac{\epsilon^2d}{3600}\right) + 2n\exp\left(-\frac{\epsilon^2d}{9000}\right).
		\end{equation}
		\item\label{it:Kemenysecond} \emph{Second order for $H^i$.} If in addition $d \geq n^{\delta}$ for some fixed $\delta\in(0,1]$, then
		\begin{equation}\label{eq:kemenysecondorder}
			H^i - (n-2) \xrightarrow{\overrightarrow{\P}_{n,p}} 0 \qquad \text{as } n\to\infty.
		\end{equation}
		\item\label{it:Hijnormalised} \emph{Normalised form for $H_{ij}$.} For every $\epsilon\in(0,1]$ and all $n \geq 300/\epsilon^2$,
		\begin{equation}\label{eq:p2phittingnormalised}
			\overrightarrow{\P}_{n,p}\left(\exists\, i\neq j: \left|\pi_jH_{ij} - 1\right| > \epsilon\right) \leq 4n\exp\left(-\frac{\epsilon^2d}{3600}\right) + 2n\exp\left(-\frac{\epsilon^2d}{9000}\right).
		\end{equation}
		\item\label{it:Hijunnormalised} \emph{Unnormalised form for $H_{ij}$.} For every $\epsilon > 0$,
		\begin{equation}\label{eq:p2phitting}
			\overrightarrow{\P}_{n,p}\left(\exists\, i\neq j: \left|\frac{1}{n}H_{ij} - 1\right| > \epsilon\right) \longrightarrow 0 \qquad\text{as } n\to\infty,
		\end{equation}
		moreover, with $\epsilon'\defeq\min\{\epsilon/4,1/4\}$ and for all $n \geq 300/\epsilon'^2$ the left-hand side is bounded by
		\begin{equation}\label{eq:p2phittingbound}
			4n\exp\left(-\frac{\epsilon'^2d}{3600}\right) + 2n\exp\left(-\frac{\epsilon'^2d}{9000}\right) + \overrightarrow{\P}_{n,p}\left(\cU_{\epsilon'}^\complement\right),
		\end{equation}
		where $\cU_{\epsilon'} \defeq \left\{\max_{j\in[n]}\left|n\pi_j - 1\right| < \epsilon'\right\}$ is the event of uniform stationary concentration of Lemma~\ref{lem:statlinfty}.
	\end{enumerate}
	The right-hand sides in parts~\ref{it:Hitarget}, \ref{it:Hijnormalised} and~\ref{it:Hijunnormalised} converge to zero as $n \to \infty$.
\end{theorem}

Since $H^i$ does not depend on $i$, parts~\ref{it:Hitarget}
and~\ref{it:Kemenysecond} are laws of large numbers for Kemeny's constant. The
constant $-2$ is read off from a trace identity. With $R \defeq P - \mathbf{1}_n\pi^t$,
on the event $\Vert R\Vert<1$,
\begin{equation}\label{eq:kemenyintro}
	H^i = \operatorname{Tr}\left(I - \mathbf{1}_n\pi^t\right) + \operatorname{Tr}(R) + \sum_{m\geq2}\operatorname{Tr}\left(R^m\right),
\end{equation}
where the first trace is $n-1$, the second is $\operatorname{Tr}(P)-1 = -1$ because
the graph has no loops, and for $d\geq n^\delta$ the sum tends to zero in
probability (Section~\ref{sec:kemenysecond}).

Part~\ref{it:Hijnormalised} is the natural point-to-point statement. Its bound is
explicit and uses only degree concentration and a norm bound for the centred
adjacency matrix. Part~\ref{it:Hijunnormalised} adds the uniform concentration of
$\pi$ due to Cooper and Frieze~\cite{cooper_frieze_2012}, which is the only input
without an explicit rate. By Remark~\ref{rem:normalisedvsnot}, the deviation of
$H_{ij}/n$ from $1$ is asymptotically that of $n\pi_j$. Parts~\ref{it:Hijnormalised}
and~\ref{it:Hijunnormalised} transfer to the average starting hitting time. Since
$H_{jj}=0$,
\begin{equation}\label{eq:Hjfromp2p}
	\pi_jH_j - 1 = \sum_{i\neq j}\pi_i\left(\pi_jH_{ij}-1\right) - \pi_j,
	\qquad
	\frac{H_j}{n} - 1 = \sum_{i\neq j}\pi_i\left(\frac{H_{ij}}{n}-1\right) - \pi_j,
\end{equation}
so $|\pi_jH_j-1|\leq\max_{i\neq j}|\pi_jH_{ij}-1|+\pi_j$, and similarly for
$H_j/n$. On the events of the proofs $\max_j\pi_j\leq 2/\sqrt n$. Hence for every
$\epsilon>2/\sqrt n$ the probability that $\max_j|\pi_jH_j-1|>\epsilon$ is at most
the right side of~\eqref{eq:p2phittingnormalised} with $\epsilon$ replaced by
$\epsilon-2/\sqrt n$. The same holds for $\max_j|H_j/n-1|$
and~\eqref{eq:p2phittingbound}, with $\epsilon$ replaced by $\epsilon-2/\sqrt n$
also in $\epsilon'$.

The hypothesis $\log(n)\ll d$ places the graph well above the threshold
$d = \log(n)$ for strong connectivity. For $d\geq C\log(n)$ with $C$ fixed, the
bounds of parts~\ref{it:Hitarget} and~\ref{it:Hijnormalised} apply to every
$\epsilon>\sqrt{9000/C}$ (Remark~\ref{rem:largeconstantregime}).

The proofs of parts~\ref{it:Hitarget} and~\ref{it:Hijnormalised} use the group
inverse of $L$, written as a Neumann series in $P-\mathbf{1}_n\pi^t$, whose norm
is less than $1$ on the good events. The argument needs $P-\mathbf{1}_n\pi^t$ to be
a contraction. This can fail for a chain that rarely moves between blocks of
states, where the singular value formula still applies: it is evaluated in closed
form on lumped chains with $\Vert P - \mathbf{1}_n\pi^t\Vert>1$
(Example~\ref{ex:lumped} and Remark~\ref{rem:contractionfails}). Part~\ref{it:Hitarget} also follows from
Theorem~\ref{thm:Hittingsingular}, through the location of the singular values of
$L$ and the bound $(n-1)-\sum_kv_k^tu_k\leq2\Vert P-\frac1nJ\Vert_F^2$
(Appendix~\ref{sec:svdroute}). Part~\ref{it:Kemenysecond} uses the traces of
the powers of $P$.

The singular value formula involves all singular values of $L$. The next theorem
gives the limit of the empirical distribution of their centred and rescaled
values.

\begin{theorem}[Semicircle law for singular values]
	\label{thm:semicircle}
	Consider a simple random walk on a directed Erd\H{o}s--R\'enyi random graph $\overrightarrow{G}(n,p_n)$ and assume $\log(n) \ll d$ and $\limsup_{n}p_n<1$. Let $L = I-P$ be the generator of the walk, and let
	$0 = \sigma_n(L) \leq \sigma_{n-1}(L) \leq \dots \leq \sigma_1(L)$
	be its singular values. The empirical spectral distribution (ESD) of the centred and scaled singular values, defined by
	\begin{align}
		\mu_n \defeq \frac{1}{n} \sum_{k=1}^{n} \delta_{\sqrt{2d/(1-p_n)}\,(\sigma_k(L) - 1)},
	\end{align}
	converges weakly in probability to the Wigner semicircle distribution $\mu_{sc}$, whose density is given by $\rho_{sc}(x) = \frac{1}{2\pi}\sqrt{4-x^2}$ for $x \in [-2, 2]$ and $\rho_{sc}(x)=0$ otherwise.
\end{theorem}

If $p_n\to0$, the scaling is $\sqrt{2d}$. The factor $(1-p_n)^{-1/2}$ normalises
the entry variance of the adjacency matrix, as in the semicircle law for the
normalised Laplacian of the undirected graph~\cite{chen2025semicircle}. The
hypothesis $\log(n)\ll d$ gives uniform degree concentration and makes the matrix
Bernstein estimate effective. The proof (Section~\ref{sec:semicirclelaw}) starts
from the identity
\begin{equation}\label{eq:LtLintro}
	L^tL - I = -\left(P+P^t\right) + P^tP.
\end{equation}
Near the identity, the singular values of $L$ therefore follow the eigenvalues of
the Hermitian part $\frac12(P+P^t)$ of $P$, which is its symmetric part since $P$
is real (Remark~\ref{rem:hermitianpart}). Write $W$ for the centred
adjacency matrix, so that $P$ is close to $J/n + W/d$. The leading fluctuation of
$L^tL-I$ is $-(W+W^t)/d$. Its off-diagonal entries form a Wigner array and its
deterministic diagonal is negligible. The mean $J/n$ contributes a matrix of
rank at most two, and the remaining terms are small in normalised Frobenius norm.
The eigenvalues of the adjacency matrix of $\overrightarrow{G}(n,p_n)$ are complex
and follow a circular law~\cite{basak2019circular,he2023edge}. Symmetrisation
doubles the entry variance, and this is the factor $\sqrt2$ in the scaling. A
heuristic comparison with a normal generator of the same spectrum gives a
semicircle of radius $\sqrt2$ instead of $2$ (Remark~\ref{rem:hermitianpart}).

\subsection{Related work}

The reversible spectral formula and its consequences are surveyed
in~\cite{lovasz_random_1996, aldous_reversible_2002, levin_markov_2017}.
Hitting, commute and cover times of random walks on Erd\H{o}s--R\'enyi graphs are
treated in~\cite{lowe_hitting_2014, helali2019hitting}, through effective
resistance in~\cite{sylvester2021random}, and multiple walks on random regular
graphs in~\cite{cooper_multiple_2010}. All of these walks are reversible. For
directed Erd\H{o}s--R\'enyi networks, Tishby, Biham and
Katzav~\cite{tishby2017directedhitting} study the distribution of first hitting
times. Von Luxburg, Radl and Hein~\cite{luxburg2014hitting} show
$H_{ij} = \mathrm{vol}(G)/d_j\,(1+o(1))$ uniformly over pairs on large random
neighbourhood graphs, that is $\pi_jH_{ij}\to1$, through the reversible spectral
representation. Part~\ref{it:Hijnormalised} is the directed counterpart.
Chatterjee~\cite{chatterjee_spectral_2023} bounds mixing quantities by a singular
value. We give an exact identity for a hitting time in terms of the full singular
system.

Semicircle laws for adjacency matrices of sparse undirected Erd\H{o}s--R\'enyi
graphs are classical, with local laws and edge
universality~\cite{he2019local, alt_extremal_2021}. Chen, Chen and
Zhu~\cite{chen2025semicircle} prove a semicircle law for the normalised Laplacian
$I_n - D^{-1/2}AD^{-1/2}$ of the undirected graph $G(n,p_n)$ under $np_n\to\infty$
and $\limsup_n p_n<1$, with the Moore--Penrose convention at isolated vertices,
after centring at $I_n$
and scaling by $\sqrt{np_n/(1-p_n)}$; see
also~\cite{jiang2012empirical,chi2016random,akarapipattana2023random}. This is a statement
about eigenvalues of a symmetric matrix. Similarity does not preserve singular
values, so it does not describe the singular values of a walk generator even in
the undirected case, and for directed graphs the two matrices are not similar.
Hazra and Passuello~\cite{hazra2026directed} obtain the limiting singular value
distribution of the \emph{adjacency} matrix of directed inhomogeneous random
graphs, a free multiplicative convolution built from the Marchenko--Pastur
law~\cite{marchenko_distribution_1967}. That is a global law at scale $\sqrt d$.
Theorem~\ref{thm:semicircle} is a fluctuation statement for the generator at scale
$1/\sqrt d$ around $1$, and it is governed by the Hermitian part of the noise
rather than by its singular values. We are not aware of an earlier semicircle law
for the centred singular values of a non-normal random walk generator.

\subsection{Organisation}

Section~\ref{sec:preliminaries} fixes notation and recalls the reversible
formula. Section~\ref{sec:hittingtimessvd} proves Theorem~\ref{thm:Hittingsingular}
and gives the examples. Sections~\ref{sec:applicationdirectedER}
and~\ref{sec:semicirclelaw} prove Theorems~\ref{thm:hittingdirectedER}
and~\ref{thm:semicircle}. Section~\ref{sec:discussion01} collects remarks on the results.
Appendix~\ref{app:tools} collects standard tools, Appendix~\ref{app:proofs}
routine proofs, Appendix~\ref{sec:svdroute} the singular value proof of
part~\ref{it:Hitarget}, and Appendix~\ref{sec:circulantexample} the oriented
cycle.

\section{Preliminaries}\label{sec:preliminaries}
\subsection{Notation}\label{sec:notation}

Throughout, $[n] \defeq \{1,\dots,n\}$. For a matrix $A \in \C^{m\times n}$ we write
$\sigma_1(A) \geq \sigma_2(A) \geq \cdots$ for its singular values in non-increasing
order, $\Vert A\Vert$ for the spectral norm, $\Vert A\Vert_F$ for the Frobenius norm,
and $A^*$ for the conjugate transpose. We consider a discrete-time Markov chain
$(X_t)_{t\in\N}$ on a finite state space of size $n$ with transition matrix $P$,
stationary distribution $\pi$, and generator
\begin{equation}\label{eq:generator}
	L \defeq I - P.
\end{equation}
The similarity transformed generator is
\begin{equation}\label{eq:Lhat}
	\hat L \defeq \operatorname{diag}(\pi)^{1/2} L \operatorname{diag}(\pi)^{-1/2}.
\end{equation}
We write $\mathbf{1}_n$ for the vector of ones of length $n$ and
$J\defeq\mathbf{1}_n\mathbf{1}_n^t$. All Markov matrices and singular vectors in this paper are real, and for a real
matrix or vector $^t$ denotes the transpose. Complex vectors, such as the
eigenvectors of a non-normal $P$, use $^*$. Vectors are columns, so $\pi^tP = \pi^t$, and
$\mathbf{1}_n\pi^t$ is the rank-one matrix each of whose rows is $\pi^t$.
We write $a_n \ll b_n$ or $a_n = o(b_n)$ if $a_n/b_n\to0$, and
$a_n = \Theta(b_n)$ if $a_n/b_n$ is bounded away from $0$ and $\infty$. Asymptotic
statements are as $n\to\infty$, and an event holds \emph{with high probability} if
its probability tends to $1$.

\subsection{The directed Erd\H{o}s--R\'enyi random graph}\label{sec:directedERprelim}

\begin{definition}[Directed Erd\H{o}s--R\'enyi random graph]\label{def:directedER}
	Let $n \in \N$ and $p \in [0,1]$. The directed Erd\H{o}s--R\'enyi random graph
	$\overrightarrow{G}(n,p)$ is the random digraph on the vertex set $V = [n]$ in
	which, for each ordered pair $(i,j) \in V \times V$ with $i \neq j$, the arc
	$i \to j$ is present independently with probability $p$; there are no
	self-loops. We write $\overrightarrow{\P}_{n,p}$ and $\overrightarrow{\E}_{n,p}$
	for the corresponding probability measure and expectation.
\end{definition}

For $v \in V$ we write $d_v^+$ and $d_v^-$ for the out-degree and the in-degree.
Their mean is $(n-1)p$, and we work with
\begin{equation}\label{eq:ddef}
	d \defeq np
\end{equation}
instead. The relative error $1/n$ is absorbed into the degree concentration event
of Section~\ref{sec:applicationdirectedER}. Likewise
$\overrightarrow{\E}_{n,p}[A] = p(J-I)$, and the correction $pI$ has spectral norm
$d/n$ and squared Frobenius norm $np^2$, which are negligible below. The simple
random walk has transition matrix $P = D_+^{-1}A$, where $A$ is the adjacency
matrix and $D_+ = \operatorname{diag}(d_v^+)$. It is not reversible, and the
eigenvalues of $P$ are in general complex. For $\log(n) \ll d$ the graph is
strongly connected with high probability, so $P$ is irreducible and $\pi$ is
unique and positive. As a fallback convention, we set $P\defeq I$ if some out-degree is zero. On the
complement of the strong connectivity event we set
$\pi\defeq\frac1n\mathbf{1}_n$ and all hitting times equal to $0$. Every
statement about $\pi$ and the hitting times is proved on an event contained in
the strong connectivity event (Lemma~\ref{lem:stronglyconnected}), and every
statement about $P$ alone on an event where all out-degrees are positive, so this
convention plays no role. In the random graph analysis we work with $L$ rather than
$\hat L$: $L$ exposes the degree and adjacency fluctuations directly, whereas
$\hat L$ involves the random stationary distribution. Both give exact hitting
time formulas (Remark~\ref{rem:Ltilde}).

\subsection{Relaxation times and hitting times}\label{sec:hittingprelim}

For reversible chains the relaxation time is $(1-\lambda_2)^{-1}$, where
$\lambda_2$ is the second largest eigenvalue of $P$. For non-reversible chains
Chatterjee~\cite{chatterjee_spectral_2023} defines the
\emph{$\sigma$-relaxation time} $t_{\mathrm{rel}}^{\sigma}\defeq 1/\sigma_{n-1}(\hat L)$.
For reversible chains $\sigma_{n-i+1}(\hat L) = 1-\lambda_i(P)$, so this is the
classical relaxation time. The singular values of $\hat L$ lie in $[0,2]$. Indeed,
by Jensen's inequality and $\pi^tP=\pi^t$,
$\Vert P f \Vert_{\ell^2(\pi)}^2 \leq \sum_i\pi_i\sum_jP_{ij}f_j^2 = \Vert f \Vert_{\ell^2(\pi)}^2$,
so $\Vert\operatorname{diag}(\pi)^{1/2}P\operatorname{diag}(\pi)^{-1/2}\Vert \leq 1$
and $\Vert\hat L\Vert \leq 2$. The singular values of $L$ may grow with $n$.

\begin{theorem}[Spectral decomposition of expected hitting times~\cite{lovasz_random_1996}]\label{thm:lovaszhittingtimes}
	Let $P$ be the transition matrix of an irreducible reversible Markov chain. Denote by $(\lambda_{k})_{k\in[n]}$ the eigenvalues of $P$ ordered as $1=\lambda_1\geq\lambda_2\geq\cdots\geq\lambda_n$, and by $(\phi_k)_{k\in[n]}$ an orthonormal basis of eigenvectors of the symmetric matrix $\operatorname{diag}(\pi)^{1/2}P\operatorname{diag}(\pi)^{-1/2}$. Then, the expected hitting time from vertex $i$ to vertex $j$ is given by
	\begin{equation}\label{eq:Hittingeigen}
		H_{ij} = \sum_{k=2}^{n} \frac{1}{1-\lambda_k} \left(\frac{\phi_{kj}^2}{\pi_j} - \frac{\phi_{ki}\phi_{kj}}{\sqrt{\pi_i \pi_j}}\right).
	\end{equation}
\end{theorem}

\begin{definition}[Averaged hitting times]
	\label{def:averagehitting}
	The average target hitting time from $i$ and the average starting hitting time
	to $j$ are
	\begin{equation}
		H^i \defeq \sum_{j=1}^n \pi_j H_{ij}, \qquad H_j \defeq \sum_{i=1}^n \pi_i H_{ij}.
	\end{equation}
\end{definition}

\begin{corollary}[Spectral decomposition of averaged hitting times]
	\label{cor:spectraldecompaverage}
	Under the assumptions of Theorem~\ref{thm:lovaszhittingtimes}, we have
	\begin{equation}
		H^i = \sum_{k=2}^n \frac{1}{1 - \lambda_k} \quad \text{and} \quad H_j  = \frac{1}{\pi_j} \sum_{k=2}^{n} \frac{1}{1 - \lambda_k} \phi_{kj}^2.
	\end{equation}
\end{corollary}

\begin{remark}[Eigenvalue expansions beyond reversibility]\label{rem:beyondreversible}
	If $P$ is diagonalisable, with right and left eigenvectors $r_k$, $\ell_k$
	normalised by $\ell_k^*r_m = \delta_{km}$, $r_1 = \mathbf{1}_n$, $\ell_1 = \pi$,
	then the group inverse of Lemma~\ref{lem:groupinverse} is
	$L^{\#} = \sum_{k\geq2}(1-\lambda_k)^{-1}r_k\ell_k^*$, and
	$H_{ij} = (L^{\#}_{jj} - L^{\#}_{ij})/\pi_j$ gives
	\begin{equation}\label{eq:Hittingeigendiag}
		H_{ij} = \frac{1}{\pi_j}\sum_{k=2}^{n}
		\frac{\left(r_{kj} - r_{ki}\right)\overline{\ell_{kj}}}{1-\lambda_k};
	\end{equation}
	see~\cite[Eq.~(22)]{riascos_nonlocal_2020}. If
	$\hat P \defeq I-\hat L$ is normal, an orthonormal eigenbasis of $\hat P$ turns
	this into~\eqref{eq:Hittingeigen}, with moduli and complex conjugates. Without
	normality the summands carry the condition number of the eigenvector matrix.
	Perturbing the defective chain of Example~\ref{ex:defective} by $\epsilon$, they
	grow like $\epsilon^{-1/2}$ while $H_{ij}$ stays bounded. Without
	diagonalisability there is no such expansion. Kemeny's constant needs neither
	hypothesis. The eigenvalues of $L^{\#}$ are $0$ and $(1-\lambda_k)^{-1}$,
	$k\geq2$, with algebraic multiplicities, so $H^i = \operatorname{Tr}(L^{\#})$ gives
	\begin{equation}\label{eq:Kemenygeneral}
		H^i = \sum_{k=2}^{n}\frac{1}{1-\lambda_k}
	\end{equation}
	for every irreducible chain.
\end{remark}

\begin{example}[Oriented cycle]
	\label{ex:MCtotheright}
	Let $P$ be the chain on $\Z/n\Z$ that moves one step to the right. Its
	eigenvalues $\exp(2\pi\mathrm{i}(k-1)/n)$ lie on the unit circle, so there is no
	spectral gap. Here $\pi$ is uniform, $\hat L = L$, and $LL^t$ is circulant with
	coefficients $c_0 = 2$, $c_1 = c_{n-1} = -1$. Hence
	\begin{equation}\label{eq:circulantlletevals}
		\sigma_k(L)^2 \in \left\{4 \sin^2\left(\frac{\pi}{n} (k-1)\right) : k \in [n]\right\},
	\end{equation}
	and
	\begin{equation}\label{eq:circulantsigmarelax}
		t_{\mathrm{rel}}^{\sigma} = \frac{1}{2\sin\left(\frac{\pi}{n}\right)} = \frac{n}{2\pi}(1+o(1)).
	\end{equation}
	The singular triplets are computed in Appendix~\ref{sec:circulantexample}.
\end{example}

\section{Expected hitting times via the SVD}\label{sec:hittingtimessvd}
In this section we prove Theorem~\ref{thm:Hittingsingular}. The proof follows
that of Theorem~\ref{thm:lovaszhittingtimes}. A one-step computation gives a
singular linear system for $H$. The Moore--Penrose inverse of $L$ gives one
solution, and the condition $H_{ii}=0$ selects the hitting times. Irreducibility
enters only through $\pi_i>0$ for all $i$ and $\operatorname{rank}(L) = n-1$,
that is $\sigma_{n-1}(L)>0$. Both follow from the Perron--Frobenius theorem.

\begin{proof}[Proof of Theorem~\ref{thm:Hittingsingular}]
	For $i\neq j$, conditioning on the first step gives
	$H_{ij} = 1 + \sum_{k} P_{ik}H_{kj}$. Hence $F \defeq J + PH - H$ is a
	diagonal matrix, and
	$F^t\pi = J\pi + H^t(P-I)^t\pi = \mathbf{1}_n$, so $F_{ii} = 1/\pi_i$. Thus
	\begin{equation}\label{eq:hittingtimematrix}
		(I-P) H = J-F.
	\end{equation}
	Since $L = I-P$ is singular, a solution $X$ of~\eqref{eq:hittingtimematrix} is
	determined only up to $X + \mathbf{1}_na^t$, $a\in\R^n$, and $H$ is the unique
	solution with
	\begin{equation}\label{eq:diagonalhittingtime}
		H_{ii} = 0 \text{ for all } i\in[n].
	\end{equation}
	Let $L = \sum_{k=1}^n \sigma_k u_k v_k^t$ be a singular value decomposition, so
	that the Moore--Penrose inverse of $L$ is (see~\cite{ben2003generalized})
	\begin{equation}\label{eq:Zdefn}
		L^{+} = \sum_{k=1}^{n-1} \frac{1}{\sigma_k} v_ku_k^t.
	\end{equation}
	By orthonormality, $LL^{+} = \sum_{k=1}^{n-1} u_ku_k^t = I - u_nu_n^t$. The left
	kernel of $L$ is spanned by $\pi$, so $u_n = \pi/\Vert\pi\Vert$. Since
	$(\pi^t(J-F))_i = \sum_j\pi_j - \pi_i/\pi_i = 0$ for all $i$, we get
	$LL^{+}(J-F) = J-F$. So $X \defeq L^{+}(J-F)$ solves~\eqref{eq:hittingtimematrix},
	and
	\begin{equation}
		X_{ij} = \sum_{k=1}^{n-1} \frac{1}{\sigma_k} v_{ki}\left(u_k^t\mathbf{1}_n - \frac{u_{kj}}{\pi_j}\right).
	\end{equation}
	By~\eqref{eq:diagonalhittingtime}, $a_j = -X_{jj}$, and $H_{ij} = X_{ij} - X_{jj}$
	is~\eqref{eq:Hittingsingular}.
\end{proof}

Averaging~\eqref{eq:Hittingsingular} against $\pi$ in $j$ or in $i$, and using
$\sum_j\pi_j = 1$, gives the following.

\begin{corollary}[Singular decomposition of the average target and starting hitting times]
	\label{cor:averagedhittingtimes}
	Let $P$ be the transition matrix of an irreducible Markov chain on the state space $[n]$, with stationary distribution $\pi$. With the notation $(\sigma_k)_{k\in[n]}$, $(u_k)_{k\in[n]}$ and $(v_k)_{k\in[n]}$ of Theorem~\ref{thm:Hittingsingular}, the average target hitting time from state $i$ is given by
	\begin{equation}\label{eq:meantargethittingtime}
		H^i = \sum_{j=1}^n \pi_j H_{ij} = \sum_{k=1}^{n-1} \frac{1}{\sigma_k}\left(v_k^tu_k - v_k^t\pi u_k^t\mathbf{1}_n\right),
	\end{equation}
	whereas the average starting hitting time is given by
	\begin{equation}\label{eq:meanstartinghittingtime}
		H_j = \sum_{i=1}^n \pi_i H_{ij} = \sum_{k=1}^{n-1} \frac{1}{\sigma_{k}} \left(v_k^t\pi -v_{kj}\right)\left(u_k^t\mathbf{1}_n - \frac{u_{kj}}{\pi_j}\right).
	\end{equation}
\end{corollary}

\begin{proof}
	For~\eqref{eq:meantargethittingtime}, the terms $v_{ki}u_k^t\mathbf{1}_n$ and
	$-v_{ki}\sum_j u_{kj}$ cancel, and $\sum_j\pi_jv_{kj}u_{kj}/\pi_j = v_k^tu_k$.
	For~\eqref{eq:meanstartinghittingtime}, $\sum_i\pi_iv_{ki} = v_k^t\pi$.
\end{proof}

\begin{remark}[Moore--Penrose form]\label{rem:MPform}
	Collecting terms in~\eqref{eq:Hittingsingular} through~\eqref{eq:Zdefn} gives
	\begin{equation}\label{eq:HijviaMP}
		H_{ij} = \left(L^{+}\mathbf{1}_n\right)_i - \left(L^{+}\mathbf{1}_n\right)_j + \frac{\left(L^{+}\right)_{jj} - \left(L^{+}\right)_{ij}}{\pi_j},
		\qquad
		H^i = \operatorname{Tr}\left(L^{+}\right) - \pi^tL^{+}\mathbf{1}_n.
	\end{equation}
	The group inverse $L^{\#}$ of Section~\ref{sec:groupinverse} satisfies
	$L^{\#}\mathbf{1}_n = 0$ and $\pi^tL^{\#} = 0$, whereas $L^{+}\pi = 0$ and
	$\mathbf{1}_n^tL^{+} = 0$. The two coincide only when $\pi$ is uniform. This is
	why~\eqref{eq:HijviaMP} carries the extra terms $(L^{+}\mathbf{1}_n)_i - (L^{+}\mathbf{1}_n)_j$
	and $\pi^tL^{+}\mathbf{1}_n = \sum_{k<n}\sigma_k^{-1}(v_k^t\pi)(u_k^t\mathbf{1}_n)$
	compared with Meyer's identities~\eqref{eq:groupinversehitting}.
\end{remark}

\begin{remark}[The transformed generator]\label{rem:Ltilde}
	Let $(\hat\sigma_k,\hat u_k,\hat v_k)$ be the singular triplets of $\hat L$
	in~\eqref{eq:Lhat}, and write $\Pi^{1/2}\defeq\operatorname{diag}(\pi)^{1/2}$.
	In the proof, $L^{+}$ may be replaced by any matrix $G$ with
	$LG(J-F) = J-F$. Take $G \defeq \Pi^{-1/2}\hat L^{+}\Pi^{1/2}$. The left kernel
	of $\hat L$ is spanned by the unit vector $\sqrt\pi$, so
	$\hat L\hat L^{+} = I - \sqrt\pi\sqrt\pi^t$ and
	\begin{equation}
		LG = \Pi^{-1/2}\hat L\hat L^{+}\Pi^{1/2} = I - \mathbf{1}_n\pi^t.
	\end{equation}
	Since $\pi^t(J-F) = 0$, this gives $LG(J-F) = J-F$, and the proof gives
	\begin{equation}\label{eq:hittingsingulartilde}
		H_{ij} = \sum_{k=1}^{n-1} \frac{1}{\hat{\sigma}_k} \left(\frac{\hat{v}_{ki}}{\sqrt{\pi_i}} - \frac{\hat{v}_{kj}}{\sqrt{\pi_j}}\right)\left(\hat{u}_k^t\sqrt{\pi} - \frac{\hat{u}_{kj}}{\sqrt{\pi_j}}\right).
	\end{equation}
	Both the left and the right kernel of $\hat L$ are spanned by $\sqrt\pi$, so
	$\hat u_k^t\sqrt\pi = \hat v_k^t\sqrt\pi = 0$ for $k<n$, for every chain. Hence
	\begin{equation}\label{eq:meantargethittingtimetilde}
		H^i = \sum_{k=1}^{n-1} \frac{\hat{v}_k^t\hat{u}_k}{\hat{\sigma}_k},
		\qquad
		H_j = \frac{1}{\pi_j}\sum_{k=1}^{n-1} \frac{\hat{u}_{kj}\hat{v}_{kj}}{\hat{\sigma}_{k}}.
	\end{equation}
	For reversible chains $\hat u_k = \hat v_k = \phi_{n-k+1}$ and
	$\hat\sigma_{n-k+1} = 1-\lambda_k$, so~\eqref{eq:hittingsingulartilde}
	and~\eqref{eq:meantargethittingtimetilde} reduce to
	Theorem~\ref{thm:lovaszhittingtimes} and
	Corollary~\ref{cor:spectraldecompaverage}. For the directed
	Erd\H{o}s--R\'enyi graph we work with $L$, on which the degree and adjacency
	estimates are stated.
\end{remark}

\begin{remark}[Orientation changes the order of magnitude]\label{rem:cyclecomparison}
	For the oriented cycle of Example~\ref{ex:MCtotheright},
	Corollary~\ref{cor:averagedhittingtimes} gives $H^i = H_j = (n-1)/2$
	(Appendix~\ref{sec:circulantexample}). Here $P$ is normal and the eigenvalue sum
	$\sum_{k=1}^{n-1}(1-e^{2\pi\mathrm{i}k/n})^{-1} = (n-1)/2$ gives the same value. For the
	simple random walk on the unoriented cycle,
	$H_{ij} = \operatorname{dist}(i,j)(n-\operatorname{dist}(i,j))$, so
	\begin{equation}\label{eq:undirectedcycleaverage}
		H^i = \frac{1}{n}\sum_{j=1}^n \operatorname{dist}(i,j)\left(n - \operatorname{dist}(i,j)\right) = \frac{n^2-1}{6}.
	\end{equation}
	For these two walks the ratio of the unoriented to the oriented value of $H^i$
	is $(n+1)/3$.
\end{remark}

\begin{example}[A defective chain]\label{ex:defective}
	Consider the directed triangle $1 \to 2 \to 3 \to 1$ with holding probability
	$3/4$ at states $2$ and $3$:
	\begin{equation}\label{eq:defectiveP}
		P = \begin{pmatrix}
			0 & 1 & 0 \\[2pt]
			0 & \tfrac{3}{4} & \tfrac{1}{4} \\[2pt]
			\tfrac{1}{4} & 0 & \tfrac{3}{4}
		\end{pmatrix},
		\qquad
		\pi = \tfrac{1}{9}\left(1, 4, 4\right).
	\end{equation}
	The chain is irreducible and aperiodic. Its characteristic polynomial is
	$(\lambda - 1)(\lambda - \tfrac14)^2$, and the eigenvalue $\tfrac14$ has
	geometric multiplicity one. So $P$ has no basis of eigenvectors, over $\R$ or
	over $\C$. The singular values of $L$ are $\sqrt{33}/4$, $\sqrt3/4$ and $0$,
	and~\eqref{eq:Hittingsingular} returns
	\begin{equation}
		\left(H_{ij}\right)_{i,j\in[3]} = \begin{pmatrix}
			0 & 1 & 5\\
			8 & 0 & 4\\
			4 & 5 & 0
		\end{pmatrix},
	\end{equation}
	These values can be read off from the walk: it moves along $1\to2\to3\to1$ and
	stays at $2$ and at $3$ for geometric times with mean $4$.
\end{example}

\subsection{Slow singular modes}\label{sec:slowmodes}

Since the left and right singular vectors form orthonormal bases, the formula can
be truncated to its slow modes with an error that involves no eigenvector
condition number.

\begin{proposition}[Truncation to slow singular modes]\label{prop:slowmodes}
	Let $P$ be irreducible, let $1\leq s<n-1$, and write $S \defeq \{n-s,\dots,n-1\}$ for the indices of the $s$ smallest non-zero singular values of $L$. Then for all $i,j\in[n]$,
	\begin{equation}\label{eq:slowmodes}
		\left|H_{ij} - \sum_{k\in S}\frac{1}{\sigma_k}\left(v_{ki}-v_{kj}\right)\left(u_k^t\mathbf{1}_n - \frac{u_{kj}}{\pi_j}\right)\right| \leq \frac{\sqrt2}{\sigma_{n-s-1}}\left(\sqrt n + \frac{1}{\pi_j}\right).
	\end{equation}
\end{proposition}
\begin{proof}
	Write $w \defeq \mathbf{1}_n - e_j/\pi_j$, so that the summand in~\eqref{eq:Hittingsingular} is $\sigma_k^{-1}\left(v_k^t(e_i-e_j)\right)\left(u_k^tw\right)$. The modes $k\leq n-s-1$ have $\sigma_k\geq\sigma_{n-s-1}$, and by Cauchy--Schwarz their total contribution is at most $\sigma_{n-s-1}^{-1}\Vert V^t(e_i-e_j)\Vert\,\Vert U^tw\Vert \leq \sigma_{n-s-1}^{-1}\sqrt2\,\Vert w\Vert$, with $\Vert w\Vert\leq\sqrt n + 1/\pi_j$.
\end{proof}

The bound is of order $1/(\pi_j\sigma_{n-s-1})$, which is the order of the fast
part on lumped chains (Example~\ref{ex:lumped}). The truncation is therefore
informative when the slow modes carry hitting times much larger than $1/\pi_j$.
On the directed Erd\H{o}s--R\'enyi graph all non-zero singular values are of order
one, and no truncation is needed.

An eigenvalue truncation of the same kind carries the condition number of the
eigenvector matrix, which is unbounded over non-normal chains and infinite for
defective ones. In the next example the slow and fast contributions are computed
in closed form.

\begin{example}[Lumped chains]\label{ex:lumped}
	Let $M$ be an irreducible stochastic $r\times r$ matrix with stationary
	distribution $\pi_M$ and hitting times $H^M$, and let
	$P \defeq M\otimes\frac1mJ_m$ on $[r]\times[m]$. From $(A,x)$ the walk moves to
	block $B$ with probability $M_{AB}$ and then to a uniform state of that block.
	Then $\pi = \pi_M\otimes\frac1m\mathbf{1}_m$. The space $\R^n$ splits
	orthogonally into $\{a\otimes\mathbf{1}_m\}$ and
	$F \defeq \bigoplus_A e_A\otimes\mathbf{1}_m^\perp$, and both are invariant under
	$L$ and $L^t$. On the first, $L$ acts as $I-M$, so each singular triplet
	$(\tilde\sigma,\tilde v,\tilde u)$ of $I-M$ gives the triplet
	$(\tilde\sigma,\tilde v\otimes\mathbf{1}_m/\sqrt m,\tilde u\otimes\mathbf{1}_m/\sqrt m)$
	of $L$. On $F$, $L$ is the identity, which gives $r(m-1)$ singular values equal
	to $1$ with $u_k = v_k$ and $u_k^t\mathbf{1}_n = 0$. We choose the singular
	vectors of $L$ in this way, adapted to the two invariant subspaces; this matters
	only if $I-M$ has a singular value equal to $1$. In~\eqref{eq:Hittingsingular}
	the $r-1$ block modes with non-zero singular value reproduce
	Theorem~\ref{thm:Hittingsingular} for $M$, and the modes in $F$ contribute
	$((\Pi_F)_{jj} - (\Pi_F)_{ij})/\pi_j$, with $\Pi_F$ the orthogonal projection
	onto $F$. Since $\Pi_F = I_r\otimes(I_m-\frac1mJ_m)$ and
	$\pi_{(B,y)} = \pi_M(B)/m$, this gives
	\begin{equation}\label{eq:lumpedhitting}
		H_{(A,x),(B,y)} = H^M_{AB} + \frac{m-1}{\pi_M(B)} \quad (A\neq B), \qquad H_{(A,x),(A,y)} = \frac{m}{\pi_M(A)} \quad (x\neq y).
	\end{equation}
	The block contribution is $H^M_{AB}$, and the contribution of $F$ is of order
	$1/\pi_j$. If all non-zero singular values of $I-M$ lie strictly below $1$, these
	are the slow and fast contributions of Proposition~\ref{prop:slowmodes} with
	$s = r-1$. Nothing here requires $M$ to be reversible, normal
	or diagonalisable.
\end{example}

\section{Application to directed Erd\H{o}s--R\'enyi graphs}\label{sec:applicationdirectedER}
In this section we prove Theorem~\ref{thm:hittingdirectedER}. We introduce two
good events. The event $\cF_{\nu}$ says that all in- and out-degrees lie within
$\nu d$ of $d$, and $\cH_{\nu}$ bounds the norm of the centred adjacency matrix.
On $\cF_{\nu_1}\cap\cH_{\nu_2}$ the stationary distribution is close to uniform in
$\ell^2$, and $\Vert P - \mathbf{1}_n\pi^t\Vert$ is bounded away from $1$. The
group inverse of $L$ is then a Neumann series. Its entries give
parts~\ref{it:Hijnormalised} and~\ref{it:Hijunnormalised}, and its trace gives
part~\ref{it:Hitarget}. Part~\ref{it:Kemenysecond} controls the traces of the
series term by term, through the moments of the number of closed walks.

If $\nu_1+\nu_2<1/2$, then on $\cF_{\nu_1}\cap\cH_{\nu_2}$ the digraph is
strongly connected, so $P$ is irreducible and $\pi$ is unique and positive
(Lemma~\ref{lem:stronglyconnected} in Appendix~\ref{sec:connectivityproof}). This
is deterministic, so it adds no term to the probability bounds. The proofs below
take $\nu_1=\nu_2\leq1/30$, or $\nu_1,\nu_2\to0$, so the condition holds
throughout.

\subsection{Degree concentration and adjacency norm bound}

Let $\cF^+_{\nu}$ and $\cF^-_{\nu}$ be the events
\begin{equation}
	\cF^+_{\nu} \defeq \left\{\forall k\in V, \left|\frac{1}{d}\sum_{\ell \in V} A_{k\ell} - 1\right| \leq \nu\right\},
	\qquad
	\cF^-_{\nu} \defeq \left\{\forall k\in V, \left|\frac{1}{d}\sum_{\ell \in V} A_{\ell k} - 1\right| \leq \nu\right\},
\end{equation}
and set
\begin{equation}\label{eq:Fevent}
	\cF_{\nu} \defeq \cF^-_{\nu} \cap \cF^+_\nu.
\end{equation}

\begin{lemma}[Concentration of in- and out-degrees]
	\label{lem:degreeconcentrationinout}
	For any $\nu \in (0,1)$ and any $n \geq 10/\nu$,
	\begin{equation}\label{eq:degreeconcentrationout}
		\overrightarrow{\P}_{n,p}\left(\left(\cF^+_{\nu}\right)^\complement\right) \leq 2n\exp\left(-\frac{\nu^2d}{4}\right)
	\end{equation}
	and
	\begin{equation}\label{eq:degreeconcentrationin}
		\overrightarrow{\P}_{n,p}\left(\left(\cF^-_{\nu}\right)^\complement\right) \leq 2n\exp\left(-\frac{\nu^2d}{4}\right).
	\end{equation}
\end{lemma}

The proof is a Chernoff bound and a union bound (Appendix~\ref{sec:concentrationproofs}).
The constant $4$ instead of $3$ absorbs the difference between the mean degree
$(n-1)p$ and $d$.

\begin{corollary}\label{cor:degreeconcentration}
	For any $\nu \in (0,1)$ and any $n\geq 10/\nu$,
	\begin{equation}
		\overrightarrow{\P}_{n,p}\left(\cF_{\nu}^\complement\right) \leq 4n\exp\left(-\frac{\nu^2d}{4}\right).
	\end{equation}
\end{corollary}

Let
\begin{equation}\label{eq:Wdef}
	W \defeq A - pJ.
\end{equation}
Since $\E(A) = p(J-I)$, we have $W = (A - \E A) - pI$. We centre by $pJ$ because
$W/d = A/d - J/n$ is the combination that occurs below. For $\nu > 0$ let
\begin{equation}\label{eq:Hnuevent}
	\cH_{\nu} \defeq \left\{\Vert W\Vert \leq \nu d\right\}.
\end{equation}

\begin{lemma}[High probability bound on the norm of $W$]
	\label{lem:Wignernormbound}
	For any $\nu > 0$ and all $n \geq 2/\nu$, we have
	\begin{equation}\label{eq:Wignernorm}
		\overrightarrow{\P}_{n,p}\left(\cH_{\nu}^\complement\right) \leq 2n\exp\left(-\frac{\nu^2 d}{8 + \frac{4\nu}{3}}\right).
	\end{equation}
\end{lemma}

The proof applies the matrix Bernstein inequality (Appendix~\ref{sec:concentrationproofs}).
For $0 < \nu < 3/2$ the bound becomes
\begin{equation}\label{eq:Wignernorm2}
	\overrightarrow{\P}_{n,p}\left(\cH_{\nu}^\complement\right) \leq 2n\exp\left(-\frac{\nu^2 d}{10}\right).
\end{equation}

\begin{corollary}[From good events to high probability]\label{cor:goodprobability}
	Let $\nu_1\in(0,1)$ and $0 < \nu_2 < 3/2$, and let $A$ be an event with $\cF_{\nu_1}\cap\cH_{\nu_2} \subseteq A$. Then, for all $n \geq \max\left\{10/\nu_1, 2/\nu_2\right\}$,
	\begin{equation}\label{eq:goodprobability}
		\overrightarrow{\P}_{n,p}\left(A^\complement\right) \leq 4n\exp\left(-\frac{\nu_1^2d}{4}\right) + 2n\exp\left(-\frac{\nu_2^2d}{10}\right).
	\end{equation}
\end{corollary}

\begin{remark}[Fixed multiples of $\log n$]\label{rem:largeconstantregime}
	The proofs below take $\nu_1 = \nu_2 = \epsilon/30$, and then the right-hand side
	of~\eqref{eq:goodprobability} tends to zero if
	\begin{equation}\label{eq:Cthreshold}
		d \geq C\log(n) \qquad\text{with}\qquad C > \frac{9000}{\epsilon^2}.
	\end{equation}
	For a fixed $C$ the bounds of parts~\ref{it:Hitarget} and~\ref{it:Hijnormalised}
	therefore apply to every $\epsilon>\sqrt{9000/C}$, with polynomially small
	right-hand sides. Every $\epsilon>0$ is covered under $\log(n)\ll d$.
\end{remark}

\subsection{Concentration of the stationary distribution}

Write
\begin{equation}\label{eq:xdefstat}
	x \defeq \pi - \frac{1}{n}\mathbf{1}_n,
\end{equation}
so that $x^t\mathbf{1}_n = 0$ and $n\Vert\pi\Vert^2 = 1 + n\Vert x\Vert^2$. The
$\ell^2$ estimate below is deterministic on the good events. The uniform estimate
is finer. An $\ell^2$ bound of order $n^{-1/2}$ does not give the coordinatewise
order $n^{-1}$, and for that we use the analysis of Cooper and
Frieze~\cite{cooper_frieze_2012}.

\begin{lemma}[Concentration of the stationary distribution in $\ell^2$]
	\label{lem:statalignment}
	For every $\epsilon > 0$, there exist constants $\nu_1, \nu_2 > 0$ such that, on the event $\cF_{\nu_1} \cap \cH_{\nu_2}$,
	\begin{equation}\label{eq:statalignmentl2}
		n\Vert x\Vert^2 < \epsilon.
	\end{equation}
\end{lemma}

\begin{proof}
	Let $Q \defeq P - J/n$. Since $x^tJ = 0$, the equation $\pi^tP = \pi^t$ reads
	$\frac1n\mathbf{1}_n^tQ + x^tQ = x^t$, so that, whenever $\Vert Q\Vert<1$,
	\begin{equation}\label{eq:xexpression}
		x^t = \frac{1}{n}\mathbf{1}_n^tQ\left(I - Q\right)^{-1}.
	\end{equation}
	We have $\Vert Q\Vert \leq \Vert P - A/d\Vert + \Vert W\Vert/d$, and
	$\Vert P - A/d\Vert^2 \leq \Vert P - A/d\Vert_1\Vert P - A/d\Vert_\infty$
	(maximal column sum and maximal row sum). On $\cF_{\nu_1}$,
	$|1/d_i^+ - 1/d| \leq \nu_1/((1-\nu_1)d)$, and every row and column of $A$ sums
	to at most $(1+\nu_1)d$. Hence both factors are at most
	$\nu_1(1+\nu_1)/(1-\nu_1)$, and on $\cF_{\nu_1}\cap\cH_{\nu_2}$
	\begin{equation}\label{eq:Qnormbound}
		\Vert Q\Vert \leq \eta, \qquad \eta \defeq \nu_1 \frac{1+\nu_1}{1 - \nu_1} + \nu_2.
	\end{equation}
	Choose $\nu_1,\nu_2$ with $\eta<1$. The column sums of $Q$ are
	$(\mathbf{1}_n^tQ)_j = \sum_iA_{ij}/d_i^+ - 1$, which lies between
	$(1-\nu_1)/(1+\nu_1) - 1$ and $(1+\nu_1)/(1-\nu_1) - 1$. Hence
	\begin{equation}\label{eq:colsumsQ}
		\left|\left(\mathbf{1}_n^tQ\right)_j\right| \leq \frac{2\nu_1}{1 - \nu_1},
		\qquad
		\left\Vert Q^t\mathbf{1}_n\right\Vert \leq \sqrt{n}\,\frac{2\nu_1}{1-\nu_1}.
	\end{equation}
	With~\eqref{eq:xexpression} and $\Vert(I-Q)^{-1}\Vert\leq(1-\eta)^{-1}$,
	\begin{equation}\label{eq:xl2bound}
		\Vert x\Vert \leq \frac{1}{\sqrt{n}}\,\frac{2\nu_1}{\left(1-\nu_1\right)\left(1-\eta\right)}.
	\end{equation}
	The right-hand side of $n\Vert x\Vert^2 \leq (2\nu_1/((1-\nu_1)(1-\eta)))^2$
	tends to $0$ as $\nu_1,\nu_2\to0$.
\end{proof}

\begin{lemma}[Uniform concentration of the stationary distribution]
	\label{lem:statlinfty}
	Assume $d = np_n$ satisfies $\log(n) \ll d$, and
	for $\epsilon > 0$ define the event
	\begin{equation}\label{eq:Uevent}
		\cU_{\epsilon} \defeq \bigcap_{j\in[n]}\left\{\left|n\pi_j - 1\right| < \epsilon\right\}.
	\end{equation}
	Then $\overrightarrow{\P}_{n,p}\left(\cU_{\epsilon}\right) \to 1$ as $n\to\infty$, for every $\epsilon > 0$.
\end{lemma}

\begin{proof}
	The model of~\cite{cooper_frieze_2012} is Definition~\ref{def:directedER}. Its
	expected number of arcs is $m = n(n-1)p$, which matches the loopless convention
	of Definition~\ref{def:directedER}. Loops would
	in any case change each degree by at most one, which is $o(d)$, and the
	conclusion below would be unchanged.
	Writing $\gamma \defeq np - \log n$, \cite[Theorem 2]{cooper_frieze_2012} gives,
	for $\gamma\to\infty$, with high probability and uniformly in $v \in V$,
	\begin{equation}\label{eq:CFthm2}
		\pi_v \sim \frac{d_v^- + \varsigma^*(v)}{m}, \qquad
		\varsigma^*(v) \defeq \max_{w \in N^-(v)} \frac{d_w^-}{d_w^+}, \qquad
		m \defeq n(n-1)p.
	\end{equation}
	where $a_v\sim b_v$ means $\max_{v}|a_v/b_v - 1|\to0$, and $N^-(v)$ is the set of
	in-neighbours of $v$. For $\log(n)\ll d$ we have $\gamma = (1+o(1))d\to\infty$.
	On $\cF_{1/2}$ every ratio $d_w^-/d_w^+$ is at most $3$, so
	$\varsigma^*(v)\leq3=o(d)$. So $\pi_v \sim d_v^-/m$ uniformly. By
	Lemma~\ref{lem:degreeconcentrationinout} with
	$\nu=\nu_n\defeq\sqrt{6\log(n)/d}$, which tends to $0$ because $\log(n)\ll d$,
	we have $d_v^- \sim d$ uniformly, and
	$m\sim nd$. Hence $n\pi_v\to1$ uniformly in probability;
	see also~\cite[Remark 3]{cooper_frieze_2012}.
\end{proof}

Lemma~\ref{lem:statlinfty} is the only ingredient that is not deterministic on
$\cF_{\nu_1}\cap\cH_{\nu_2}$. It is used only for part~\ref{it:Hijunnormalised},
which is why that part has no explicit rate.

\subsection{The group inverse and the \texorpdfstring{$\pi$}{pi}-normalised hitting times}
\label{sec:groupinverse}

For irreducible $P$ the matrix $I - P + \mathbf{1}_n\pi^t$ is invertible. Its
inverse is the fundamental matrix of Kemeny and Snell~\cite{kemeny_snell_1976},
and the group inverse of $L$ is
\begin{equation}\label{eq:groupinverse}
	L^{\#} \defeq \left(I - P + \mathbf{1}_n\pi^t\right)^{-1} - \mathbf{1}_n\pi^t.
\end{equation}

\begin{lemma}[Group inverse identities~\cite{meyer_group_1975}]\label{lem:groupinverse}
	Let $P$ be irreducible with stationary distribution $\pi$, and let $L^{\#}$ be as in~\eqref{eq:groupinverse}. Then $L^{\#}\mathbf{1}_n = 0$, $\pi^tL^{\#} = 0$, and for all $i,j\in[n]$,
	\begin{equation}\label{eq:groupinversehitting}
		H_{ij} = \frac{L^{\#}_{jj} - L^{\#}_{ij}}{\pi_j}, \qquad \pi_jH_j = L^{\#}_{jj}, \qquad H^i = \operatorname{Tr}\left(L^{\#}\right).
	\end{equation}
\end{lemma}

\begin{proof}
	The matrix $I - P + \mathbf{1}_n\pi^t$ fixes $\mathbf{1}_n$ and $\pi^t$ (from the
	left), hence so does its inverse, and $L^{\#}\mathbf{1}_n = 0$,
	$\pi^tL^{\#} = 0$ follow. The first identity is~\cite[Theorem 3.3]{meyer_group_1975}.
	Averaging it against $\pi$ in $i$ and using $\pi^tL^{\#}=0$ gives the second.
	Averaging against $\pi$ in $j$ and using $L^{\#}\mathbf{1}_n=0$ gives the third.
\end{proof}

\begin{lemma}[Uniform control of the $\pi$-normalised hitting times]\label{lem:pinormalised}
	Set $R \defeq P - \mathbf{1}_n\pi^t$ and
	\begin{equation}\label{eq:etaprime}
		\eta' \defeq \eta + \frac{2\nu_1}{\left(1-\nu_1\right)\left(1-\eta\right)},
	\end{equation}
	with $\eta$ as in~\eqref{eq:Qnormbound}. Fix $\nu_1 \in (0,1)$ and $\nu_2 > 0$ small enough that $\eta' < 1$ and $\nu_1+\nu_2<1/2$. Then, on the event $\cF_{\nu_1}\cap\cH_{\nu_2}$, we have $\Vert R\Vert \leq \eta'$ and
	\begin{equation}\label{eq:pinormalisedbound}
		\max_{j\in[n]}\left|\pi_jH_j - 1\right| \leq \frac{2}{\sqrt n} + \frac{\eta'}{1-\eta'}.
	\end{equation}
	The same argument controls the individual $\pi$-normalised hitting times: still on $\cF_{\nu_1}\cap\cH_{\nu_2}$,
	\begin{equation}\label{eq:p2pnormalisedbound}
		\max_{i\neq j}\left|\pi_jH_{ij} - 1\right| \leq \frac{4}{\sqrt n} + \frac{2\eta'}{1-\eta'}.
	\end{equation}
	Since $\eta'\to0$ as $\nu_1,\nu_2\to0$, the right-hand sides of~\eqref{eq:pinormalisedbound} and~\eqref{eq:p2pnormalisedbound} are smaller than any prescribed $\epsilon > 0$ once $\nu_1,\nu_2$ are chosen small enough and $n \geq 16/\epsilon^2$, respectively $n \geq 64/\epsilon^2$.
\end{lemma}

\begin{proof}
	Work on $\cF_{\nu_1}\cap\cH_{\nu_2}$. Since
	$\mathbf{1}_n\pi^t = \frac1nJ + \mathbf{1}_nx^t$, we have $R = Q - \mathbf{1}_nx^t$,
	and $\Vert\mathbf{1}_nx^t\Vert = \sqrt n\Vert x\Vert$. By~\eqref{eq:Qnormbound}
	and~\eqref{eq:xl2bound},
	\begin{equation}\label{eq:Rnormbound}
		\Vert R\Vert \leq \Vert Q\Vert + \sqrt n\Vert x\Vert \leq \eta'.
	\end{equation}
	Since $P\mathbf{1}_n\pi^t = \mathbf{1}_n\pi^tP = (\mathbf{1}_n\pi^t)^2 = \mathbf{1}_n\pi^t$,
	induction gives $R^m = P^m - \mathbf{1}_n\pi^t$ for $m\geq1$, and
	$I - P + \mathbf{1}_n\pi^t = I - R$. As $\Vert R\Vert<1$,
	\begin{equation}\label{eq:Gneumann}
		L^{\#} = \left(I-R\right)^{-1} - \mathbf{1}_n\pi^t = I - \mathbf{1}_n\pi^t + \sum_{m\geq1}R^m.
	\end{equation}
	Using $|N_{ij}| \leq \Vert N\Vert$ for every matrix $N$, we get
	\begin{equation}\label{eq:Gdiagbound}
		\left|L^{\#}_{jj} - 1\right| \leq \pi_j + \frac{\eta'}{1-\eta'},
		\qquad
		\left|L^{\#}_{ij}\right| \leq \pi_j + \frac{\eta'}{1-\eta'} \quad (i\neq j).
	\end{equation}
	By~\eqref{eq:xl2bound},
	$\pi_j \leq \Vert\pi\Vert \leq n^{-1/2} + \Vert x\Vert \leq (1+\eta')n^{-1/2} \leq 2n^{-1/2}$.
	Now $\pi_jH_j = L^{\#}_{jj}$ and $\pi_jH_{ij} - 1 = (L^{\#}_{jj}-1) - L^{\#}_{ij}$
	by Lemma~\ref{lem:groupinverse}, which gives~\eqref{eq:pinormalisedbound}
	and~\eqref{eq:p2pnormalisedbound}. The bound on $\Vert R\Vert$ covers all pairs at
	once, so no union bound over pairs is needed. The thresholds on $n$ make the terms
	$2/\sqrt n$ and $4/\sqrt n$ at most $\epsilon/2$.
\end{proof}

\begin{remark}[The contraction condition]\label{rem:contractionfails}
	The proof needs $\Vert P-\mathbf{1}_n\pi^t\Vert<1$, while
	Theorem~\ref{thm:Hittingsingular} does not. For the lumped chains of
	Example~\ref{ex:lumped},
	$P - \mathbf{1}_n\pi^t = (M - \mathbf{1}_r\pi_M^t)\otimes\frac1mJ_m$, and
	$\frac1mJ_m$ is an orthogonal projection, so
	$\Vert P - \mathbf{1}_n\pi^t\Vert = \Vert M - \mathbf{1}_r\pi_M^t\Vert$ for
	every $m$. This norm is not small when $M$ is close to $I$, that is, when the
	chain rarely changes block. Indeed $I-\mathbf{1}_r\pi_M^t$ is a non-zero
	idempotent, so its norm is at least $1$, and
	$\Vert M - \mathbf{1}_r\pi_M^t\Vert \geq 1-\Vert I-M\Vert$. The norm can
	exceed $1$. For
	\begin{equation}
		M = \begin{pmatrix} 1-a & a & 0\\ 0 & 1-b & b\\ c & 0 & 1-c\end{pmatrix}, \qquad (a,b,c) = (10^{-2},10^{-3},5\cdot10^{-2}),
	\end{equation}
	balancing the flow around the cycle gives
	\begin{equation*}
		\pi_M\propto(1/a,1/b,1/c), \qquad \pi_M = \tfrac{1}{1120}(100,1000,20).
	\end{equation*}
	The norm of $M - \mathbf{1}_r\pi_M^t$ is at least the Euclidean norm of its third
	row, which is
	\begin{equation*}
		\bigl\Vert\bigl(c-\pi_M(1),\,-\pi_M(2),\,1-c-\pi_M(3)\bigr)\bigr\Vert > 1.29.
	\end{equation*}
	So the Neumann series argument does not apply, for any $m$,
	while~\eqref{eq:lumpedhitting} gives the hitting times
	exactly.\footnote{Rounded to four decimal places,
	$\Vert M - \mathbf{1}_r\pi_M^t\Vert = 1.5437$.}
\end{remark}

\subsection{Proof of Theorem~\ref{thm:hittingdirectedER}}\label{sec:proofhittingER}

\begin{proof}[Proof of part~\ref{it:Hitarget}]
	Fix $\epsilon\in(0,1]$ and take $\nu_1 = \nu_2 = \nu \defeq \epsilon/30$. Since
	$\nu\leq1/10$,~\eqref{eq:Qnormbound} gives $\eta\leq 20\nu/9\leq 2/9$,
	and~\eqref{eq:etaprime} gives $\eta'\leq 20\nu/9 + 20\nu/7 < 6\nu = \epsilon/5$.
	Hence
	\begin{equation}\label{eq:explicitnu}
		\frac{\eta'}{1-\eta'} \leq \frac{\epsilon}{4},
		\qquad
		\frac{2\eta'}{1-\eta'} \leq \frac{\epsilon}{2}.
	\end{equation}
	On
	$\cF_{\nu_1}\cap\cH_{\nu_2}$, Lemma~\ref{lem:groupinverse}
	and~\eqref{eq:Gneumann} give, since $\operatorname{Tr}(I - \mathbf{1}_n\pi^t) = n-1$,
	\begin{equation}\label{eq:Hitraceseries}
		H^i = \operatorname{Tr}\left(L^{\#}\right) = (n-1) + \sum_{m\geq1}\operatorname{Tr}\left(R^m\right).
	\end{equation}
	Since $|\operatorname{Tr}(N)| \leq n\Vert N\Vert$,
	\begin{equation}\label{eq:neumannkemeny}
		\left|\frac{H^i}{n} - 1\right| \leq \frac{1}{n} + \frac{\eta'}{1-\eta'} \leq \frac{1}{n} + \frac{\epsilon}{2},
	\end{equation}
	which is at most $\epsilon$ for $n\geq2/\epsilon$. For
	$n \geq 300/\epsilon = \max\{10/\nu_1, 2/\nu_2\}$, Corollary~\ref{cor:goodprobability}
	gives~\eqref{eq:meanhittingdirectedER}. For $\epsilon>1$ the event is contained in
	the one for $\epsilon=1$.
\end{proof}

\begin{remark}\label{rem:neumann}
	Part~\ref{it:Hitarget} also follows from Theorem~\ref{thm:Hittingsingular}
	(Appendix~\ref{sec:svdroute}). The bound~\eqref{eq:neumannkemeny} puts the whole
	error into one operator norm of order $d^{-1/2}$. By part~\ref{it:Kemenysecond}
	the true deviation is $H^i - n\to-2$, so this loses a factor of order $n/\sqrt d$.
	The constant needs the traces $\operatorname{Tr}(R^m)$ themselves
	(Subsection~\ref{sec:kemenysecond}).
\end{remark}

\begin{proof}[Proof of parts~\ref{it:Hijnormalised} and~\ref{it:Hijunnormalised}]
	Part~\ref{it:Hijnormalised} follows from~\eqref{eq:p2pnormalisedbound},
	\eqref{eq:explicitnu} with $\nu_1=\nu_2=\epsilon/30$, and
	Corollary~\ref{cor:goodprobability}; the condition $n\geq300/\epsilon^2$ covers
	$n\geq64/\epsilon^2$ and $n\geq\max\{10/\nu_1,2/\nu_2\}$. For part~\ref{it:Hijunnormalised}, write,
	for $i\neq j$,
	\begin{equation}\label{eq:HijfromHijpi}
		\frac{H_{ij}}{n} - 1 = \frac{\pi_jH_{ij} - 1}{n\pi_j} + \frac{1 - n\pi_j}{n\pi_j}.
	\end{equation}
	On $\cU_{\epsilon'}$ we have $n\pi_j > 1-\epsilon'$, so
	\begin{equation}
		\max_{i\neq j}\left|\frac{H_{ij}}{n} - 1\right| \leq \frac{1}{1-\epsilon'}\left(\max_{i\neq j}\left|\pi_jH_{ij}-1\right| + \epsilon'\right).
	\end{equation}
	Given $\epsilon>0$, let $\epsilon' \defeq \min\{\epsilon/4,1/4\}$ and choose
	$\nu_1=\nu_2=\epsilon'/30$ as in part~\ref{it:Hijnormalised} for $\epsilon'$. On
	$\cF_{\nu_1}\cap\cH_{\nu_2}\cap\cU_{\epsilon'}$ and for $n\geq300/\epsilon'^2$, the
	right-hand side is at most $4\epsilon'\leq\epsilon$. A union bound
	gives~\eqref{eq:p2phittingbound}, and Lemma~\ref{lem:statlinfty}
	gives~\eqref{eq:p2phitting}.
\end{proof}

\begin{remark}\label{rem:normalisedvsnot}
	Let $a_j \defeq 1/(n\pi_j) - 1$. By~\eqref{eq:HijfromHijpi}, on $\cU_{\epsilon'}$,
	\begin{equation}\label{eq:unnormalisedequalsstationary}
		\left|\,\max_{i\neq j}\left|\frac{H_{ij}}{n}-1\right| - \max_{j\in[n]}\left|a_j\right|\,\right| \leq \frac{1}{1-\epsilon'}\max_{i\neq j}\left|\pi_jH_{ij}-1\right|,
	\end{equation}
	and the right-hand side tends to zero by part~\ref{it:Hijnormalised}. So the
	deviation of $H_{ij}/n$ from $1$ is asymptotically the deviation of $n\pi_j$ from
	$1$, and the natural statement is $\pi_jH_{ij}\approx1$.
\end{remark}

\subsection{Proof of part~\ref{it:Kemenysecond}}\label{sec:kemenysecond}

On the event where $\Vert R\Vert<1$, \eqref{eq:Hitraceseries} holds. The walk has
no loops, so $\operatorname{Tr}(R) = \operatorname{Tr}(P) - 1 = -1$, and
\begin{equation}\label{eq:kemenytraces}
	H^i - (n-2) = \sum_{m\geq2}\operatorname{Tr}\left(R^m\right) = \sum_{m\geq2}\left(\operatorname{Tr}\left(P^m\right) - 1\right).
\end{equation}
We bound the tail by a Frobenius norm estimate, and each of finitely many
remaining traces by the first two moments of the number of closed walks.

\begin{lemma}[Tail of the trace series]\label{lem:tracetail}
	Let $\nu_1\in(0,1/2]$ and $\nu_2>0$ be such that $\eta'<1$ in~\eqref{eq:etaprime}, and let $d\leq n$. Then on $\cF_{\nu_1}\cap\cH_{\nu_2}$, for every $K\geq2$,
	\begin{equation}\label{eq:tracetail}
		\sum_{m\geq K}\left|\operatorname{Tr}\left(R^m\right)\right| \leq \frac{12n}{d}\,\frac{(\eta')^{K-2}}{1-\eta'}.
	\end{equation}
\end{lemma}

\begin{proof}
	For $m\geq2$,
	\begin{equation}
		|\operatorname{Tr}(R^m)| \leq \Vert R\Vert_F\Vert R^{m-1}\Vert_F \leq \Vert R\Vert_F^2\Vert R\Vert^{m-2}.
	\end{equation}
	On $\cF_{\nu_1}$, $\Vert P\Vert_F^2 = \sum_i 1/d_i^+ \leq 2n/d$, and
	$\Vert\mathbf{1}_n\pi^t\Vert_F^2 = n\Vert\pi\Vert^2 \leq 4$ by the proof of
	Lemma~\ref{lem:pinormalised}. Hence $\Vert R\Vert_F^2 \leq 2(2n/d + 4) \leq 12n/d$,
	and we sum the geometric series.
\end{proof}

\begin{lemma}[Closed walks]\label{lem:closedwalks}
	Let $N_m \defeq \operatorname{Tr}(A^m)$ be the number of closed walks of length $m$ in $\overrightarrow{G}(n,p_n)$. For every fixed $m\geq2$ there is a constant $C_m$ such that, for all $n$ and all $1\leq d\leq n$,
	\begin{equation}\label{eq:closedwalkmoments}
		\left|\E N_m - d^m\right| \leq C_m\left(\frac{d^m}{n} + d^{m-1}\right), \qquad \operatorname{Var}\left(N_m\right) \leq C_m d^{2m-2}.
	\end{equation}
\end{lemma}

\begin{proof}
	A closed walk $w = (i_0,i_1,\dots,i_{m-1},i_0)$ contributes $\prod_l A_{i_li_{l+1}}$, whose mean is $p^{e(w)}$, with $e(w)$ the number of distinct arcs of $w$. Its arcs form a strongly connected digraph on the $k(w)$ distinct vertices of $w$, and a strongly connected digraph on $k\geq2$ vertices has out-degree at least one at every vertex and therefore at least $k$ arcs; hence $e(w)\geq k(w)$. The walks with $k(w)=m$ are the $n(n-1)\cdots(n-m+1)$ directed $m$-cycles, each with $e(w)=m$, and they contribute $d^m\prod_{l<m}(1-l/n)$. Here each cycle is counted once for each of
its $m$ starting points, as $\operatorname{Tr}(A^m)$ requires, so there is no
division by $m$. This product is $d^m$ up to $\binom m2d^m/n$. The walks with $k(w)=k<m$ fall into at most $m^m$ shapes, each realised by at most $n^k$ walks, and contribute at most $m^mn^kp^k = m^md^k \leq m^md^{m-1}$. This gives the bound on the mean.

	For the variance, $\operatorname{Var}(N_m) = \sum_{w,w'}\operatorname{Cov}\left(\prod_{w}A,\prod_{w'}A\right)$, and the covariance vanishes unless $w$ and $w'$ share an arc, in which case it is at most $p^{e(w\cup w')}$. The union of two closed walks sharing a vertex is strongly connected, so $e(w\cup w')\geq k(w\cup w')$, and two walks sharing an arc have at most $2m-2$ distinct vertices between them. For each $k$, the ordered pairs $(w,w')$ with $k$ distinct vertices fall into at most $(2m)^{2m}$ shapes, a number depending only on $m$, and each shape has at most $n^k$ labelled realisations. Since $p\leq1$ we have $p^{e(w\cup w')}\leq p^k$, and since $d\geq1$ we have $d^k\leq d^{2m-2}$. Hence $\operatorname{Var}(N_m)\leq\sum_{k\leq2m-2}(2m)^{2m}n^kp^k \leq C_md^{2m-2}$.
\end{proof}

\begin{proof}[Proof of part~\ref{it:Kemenysecond}]
	Let $d\geq n^\delta$, and note that $d = np_n\leq n$. Take
	\begin{equation}
		\nu_1 = \sqrt{\frac{12\log n}{d}}, \qquad \nu_2 = \sqrt{\frac{24\log n}{d}},
	\end{equation}
	which tend to zero because $\log(n)\ll d$. For $n$ large,
	Corollary~\ref{cor:goodprobability} gives
	\begin{equation}
		\overrightarrow{\P}_{n,p}\left((\cF_{\nu_1}\cap\cH_{\nu_2})^\complement\right) \leq 4n^{-2} + 2n^{-7/5}.
	\end{equation}
	By~\eqref{eq:Qnormbound} and~\eqref{eq:etaprime}, $\eta' \leq 4(\nu_1+\nu_2) \leq 40\sqrt{\log(n)/d}$ for $n$ large, so we may work on $\cF_{\nu_1}\cap\cH_{\nu_2}$ with $\eta'\leq1/2$ and~\eqref{eq:kemenytraces} in force.

	Fix the integer $K \defeq \lfloor 2/\delta\rfloor + 1 > 2/\delta$. By Lemma~\ref{lem:tracetail} and $d\geq n^\delta$,
	\begin{equation}
		\sum_{m\geq K}\left|\operatorname{Tr}\left(R^m\right)\right| \leq \frac{24n}{d}\left(\frac{1600\log n}{d}\right)^{\frac{K-2}{2}} \leq 24\cdot 1600^{K/2}\,(\log n)^{\frac K2}\,n^{1-\frac{\delta K}{2}} \longrightarrow 0.
	\end{equation}
	It remains to show that $\operatorname{Tr}(P^m)\to1$ in probability for each $m\in\{2,\dots,K-1\}$. Expanding the trace over closed walks and bounding each out-degree on $\cF_{\nu_1}$,
	\begin{equation}
		\frac{N_m}{\left((1+\nu_1)d\right)^m} \leq \operatorname{Tr}\left(P^m\right) = \sum_{w}\prod_{l=0}^{m-1}\frac{A_{i_li_{l+1}}}{d^+_{i_l}} \leq \frac{N_m}{\left((1-\nu_1)d\right)^m}.
	\end{equation}
	Since $\nu_1\to0$ and $m$ is fixed, it suffices that $N_m/d^m\to1$ in probability. This follows from Lemma~\ref{lem:closedwalks} and Chebyshev's inequality, since $d\to\infty$ and $\operatorname{Var}(N_m/d^m)\leq C_md^{-2}$. Inserting both estimates into~\eqref{eq:kemenytraces} proves~\eqref{eq:kemenysecondorder}.
\end{proof}

In~\eqref{eq:kemenytraces}, $n-1$ is the trace of $I-\mathbf{1}_n\pi^t$, the term
$\operatorname{Tr}(R) = -1$ is exact because there are no loops, and all other
traces are small. The polynomial lower bound on $d$ is used only for the tail.

\section{A semicircle law for the singular values of the generator}\label{sec:semicirclelaw}
In this section we prove Theorem~\ref{thm:semicircle}. The eigenvalues of $L^tL$
are the squares of the singular values of $L$. Write $\tau_n\defeq\sqrt{d/(2(1-p_n))}$,
so that $2\tau_n$ is the scaling in Theorem~\ref{thm:semicircle}. We compare
\begin{equation}
	C_n \defeq \tau_n\left(L^tL - I\right)
	\qquad\text{with}\qquad
	B_n \defeq -\tau_n\frac{W^t+W}{d} = -\frac{W^t+W}{\sqrt{2d(1-p_n)}}.
\end{equation}
The off-diagonal entries of $B_n$ form a Wigner array and its deterministic
diagonal is negligible, so its empirical spectral distribution converges to the
semicircle law (Theorem~\ref{thm:semicircleherensemble}). The difference
$C_n - B_n$ is a matrix of rank at most two plus a matrix that is small in
normalised Frobenius norm, and by Corollary~\ref{cor:esd-rank} the law transfers to
$C_n$. A Taylor expansion of the square root then transfers it to the singular
values. The sum in Theorem~\ref{thm:semicircle} includes $\sigma_n(L) = 0$. This
atom has mass $1/n$ and does not affect the limit.

\begin{figure}[h]
	\centering
	\includegraphics[width = 0.75\textwidth]{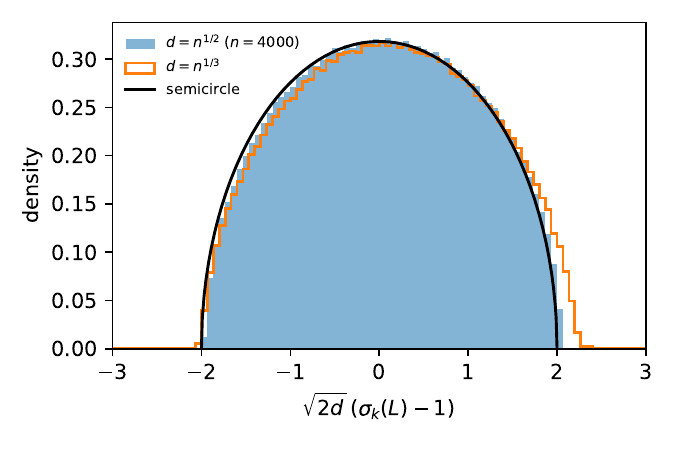}
	\caption[Semicircle histogram]{The centred, rescaled singular values $2\tau_n(\sigma_k(L)-1)$, $k<n$, of Theorem~\ref{thm:semicircle} for the simple random walk on $\overrightarrow{G}(n,p_n)$ with $n=4000$, pooled over $5$ independent realisations, for $d=n^{1/2}\approx63$ (filled) and $d=n^{1/3}\approx16$ (outline), against the semicircle density $\rho_{sc}$ (black).}
	\label{fig:semicirclehist}
\end{figure}

Figure~\ref{fig:semicirclehist} shows the law at $n=4000$. At $d=n^{1/2}$ the
histogram is close to the semicircle. At $d=n^{1/3}$, about $2\log(n)$ at this
size, it is symmetric but wider (variance about $1.06$), with soft edges. In both
cases the centre is shifted by $O(d^{-1/2})$.

\begin{remark}[Why the limit is semicircular, and why the factor $\sqrt2$]\label{rem:hermitianpart}
	For any matrix $P$,
	\begin{equation}\label{eq:hermitianpart}
		\sigma_k(I-P)^2 = \lambda_k\left(I - \left(P + P^t\right) + P^tP\right),
	\end{equation}
	so when the fluctuating part of $P$ is small the singular values of $L$ follow
	the eigenvalues of the Hermitian part $\frac12(P+P^t)$, not the singular values
	of $P$. Here the
	fluctuating part of $P$ is close to $W/d$, and the Hermitian part of a matrix with
	independent entries is a Wigner matrix. Symmetrising doubles the entrywise
	variance, whence the factor $\sqrt2$ in the scale. For an undirected graph $W = W^t$ and the
	$\sqrt2$ disappears. The following comparison is heuristic. If $L$ were normal, then $\sigma_k(L) = \vert1-\mu_k\vert$
	with $\mu_k$ the eigenvalues of $P$, and
	$\sqrt{2d}(\sigma_k(L)-1) = -\sqrt{2d}\,\operatorname{Re}\mu_k + o(1)$. The bulk
	eigenvalues follow a circular law at scale $d^{-1/2}$~\cite{bordenave2012circular},
	whose projection onto a diameter is a semicircle. This would give radius $\sqrt2$
	instead of $2$. The commutator $L^tL - LL^t = P^tP - PP^t$ is not small compared with
	the fluctuation scale $d^{-1/2}$, so $L$ is far from normal at this scale.
\end{remark}

\subsection{Expansion of the inverse degrees}

Define the diagonal matrices $\Delta$ and $\Delta_1$ by
\begin{equation}
	D_+ = d(I - \Delta), \qquad \Delta_1 \defeq \sum_{k=1}^\infty \Delta^k,
\end{equation}
so that $D_+^{-1} = d^{-1}(I + \Delta_1)$ whenever $\Vert\Delta\Vert<1$.

\begin{lemma}[Norm bounds for $\Delta_1$]
	\label{lem:normDelta1bounds}
	For any $\nu_1 \in(0,1)$, on the event $\cF_{\nu_1}$,
	\begin{equation}\label{eq:frobdelta1}
		\Vert\Delta_1\Vert \leq \frac{\nu_1}{1-\nu_1}, \qquad \Vert\Delta_1\Vert^2_F \leq n\left(\frac{\nu_1}{1-\nu_1}\right)^2.
	\end{equation}
\end{lemma}

\begin{proof}
	On $\cF_{\nu_1}$, $|\Delta_{ii}|\leq\nu_1$, so
	$|(\Delta_1)_{ii}| \leq \sum_{k\geq1}\nu_1^k = \nu_1/(1-\nu_1)$.
\end{proof}

\subsection{Convergence of the ESD of \texorpdfstring{$L^tL$}{LtL}}

\begin{lemma}[Semicircle law for $B_n$]
	\label{lem:Bsemicircle}
	Assume $d\to\infty$ and $\limsup_n p_n<1$. Then the empirical spectral distribution
	$\mu_{B_n}$ converges weakly almost surely to the Wigner semicircle law.
\end{lemma}

The proof is in Appendix~\ref{app:proofs}. The off-diagonal entries of $B_n$ have
variance exactly $1/n$, and the diagonal is the constant $2p_n/\sqrt{2d(1-p_n)}$,
which tends to zero (there are no loops).

\begin{lemma}[Semicircle law for $C_n$]
	\label{lem:eigenvaluessemicircle}
	Assume $\log(n) \ll d$ and $\limsup_n p_n<1$. Then,
	\begin{align}
		\mu_n^{\lambda} \defeq \frac{1}{n} \sum_{k=1}^{n} \delta_{\tau_n(\lambda_k(L^tL) - 1)}
	\end{align}
	converges weakly in probability to the Wigner semicircle distribution $\mu_{sc}$.
\end{lemma}

The eigenvalue $\lambda_n(L^tL)=0$ gives an atom of $\mu_n^{\lambda}$ at
$-\tau_n$, which escapes to $-\infty$. Its mass is $1/n$. Removing it changes
$\mu_n^{\lambda}$ by at most $2/n$ in the bounded Lipschitz distance, so it does
not affect the limit.

\begin{proof}
	Fix $c<1$ with $p_n\leq c$ for $n$ large, so that $\tau_n\leq\sqrt{d/(2(1-c))}$.
	By Corollary~\ref{cor:esd-rank} and Lemma~\ref{lem:Bsemicircle} it suffices
	to write $C_n-B_n = K_n + E_n$ with $\operatorname{rank}(K_n)\leq2$ and
	$\frac{1}{n}\Vert E_n\Vert_F^2 \to 0$ in probability. Let
	\begin{equation}
		\bar{\nu}_1 = \nu_1 \sqrt{\frac{\log(n)}{d}}, \qquad \bar{\nu}_2 = \nu_2 \sqrt{\frac{\log(n)}{d}},
		\qquad \nu_1 = 3,\quad \nu_2 = 4.
	\end{equation}
	Since $\log(n)\ll d\leq n$, we have $\bar\nu_1,\bar\nu_2\to0$ and
	$n\geq\max\{10/\bar\nu_1,2/\bar\nu_2\}$ for $n$ large. Corollary~\ref{cor:degreeconcentration}
	and~\eqref{eq:Wignernorm2} give
	$\overrightarrow{\P}_{n,p}\left((\cF_{\bar{\nu}_1}\cap\cH_{\bar\nu_2})^\complement\right) \leq 4n^{1-\nu_1^2/4} + 2n^{1-\nu_2^2/10}\to0$.
	We work on $\cF_{\bar{\nu}_1}\cap\cH_{\bar{\nu}_2}$.

	Insert $L = I - D_+^{-1}A$, $D_+^{-1} = d^{-1}(I+\Delta_1)$ and $A = W + pJ$ into
	$L^tL = I - A^tD_+^{-1} - D_+^{-1}A + A^tD_+^{-2}A$. With $pJ/d = J/n$ and
	$J^2 = nJ$ this gives the exact identity
	\begin{align}\label{eq:LtLexact}
		L^tL &= I - \frac{1}{n}J + \frac{2}{n^2}J\Delta_1J - \frac{1}{n}J\Delta_1 - \frac{1}{n}\Delta_1J + \frac{1}{n^2}J\Delta_1^2 J \\
		&- \frac{W^t + W + W^t\Delta_1 + \Delta_1W}{d} + \frac{JW + W^tJ + J\Delta_1^2W + W^t\Delta_1^2J + 2J\Delta_1W + 2W^t\Delta_1J}{nd} \nonumber\\
		&+\frac{1}{d^2}W^tW + \frac{2}{d^2}W^t\Delta_1W + \frac{1}{d^2}W^t\Delta_1^2W. \nonumber
	\end{align}
	Every term of~\eqref{eq:LtLexact} that contains $J$ is of the form $JX$ or
	$XJ$. The sum of the terms $JX$ has range in $\operatorname{span}\{\mathbf{1}_n\}$,
	and the sum of the terms $XJ = (X\mathbf{1}_n)\mathbf{1}_n^t$ has rank at most one.
	Let $K_n$ be $\tau_n$ times the sum of all terms containing $J$. Then
	$\operatorname{rank}(K_n)\leq 2$. The remaining terms other than
	$I$ and $-(W^t+W)/d$ give
	\begin{equation}\label{eq:Endef}
		E_n \defeq \tau_n\left(-\frac{W^t\Delta_1 + \Delta_1W}{d} + \frac{W^t(I+\Delta_1)^2W}{d^2}\right),
	\end{equation}
	and $C_n - B_n = K_n + E_n$. We use
	$\Vert A_1 \cdots A_k\Vert_F \leq \Vert A_j\Vert_F\prod_{i \neq j} \Vert A_i\Vert$,
	Lemma~\ref{lem:normDelta1bounds}, $\Vert W\Vert \leq \bar{\nu}_2 d$, and
	\begin{equation}
		\Vert W\Vert_F^2 = \sum_{i,j}\left(A_{ij}-p\right)^2
		\leq \sum_{i\neq j} A_{ij} + n^2p^2 \leq nd(1+\bar{\nu}_1) + d^2 \leq 3nd.
	\end{equation}
	For $\bar\nu_1\leq1/2$ we have $\Vert\Delta_1\Vert\leq2\bar\nu_1$ and $\Vert(I+\Delta_1)^2\Vert\leq4$. Hence
	\begin{equation}
		\Vert W^t\Delta_1+\Delta_1W\Vert_F \leq 4\bar\nu_1\sqrt{3nd},
		\qquad
		\Vert W^t(I+\Delta_1)^2W\Vert_F \leq 4\Vert W\Vert\Vert W\Vert_F \leq 4\bar\nu_2d\sqrt{3nd}.
	\end{equation}
	Multiplying by $\tau_n/d$ and by $\tau_n/d^2$ respectively, and using
	$\bar\nu_i = \nu_i\sqrt{\log(n)/d}$,
	\begin{equation}\label{eq:CBfinalbound}
		\frac{1}{n} \left\Vert E_n\right\Vert_F^2 \leq \frac{C}{1-c}\,\frac{\log(n)}{d} \longrightarrow 0
	\end{equation}
	for a numerical constant $C$.
\end{proof}

The matrix $K_n$ contains the mean $J/n$ of $P$. In normalised Frobenius norm it is
of order $\sqrt{d/n}$, which is why it is handled by the rank inequality and not by
the Frobenius bound. The factor $\sqrt{1-p_n}$ in the scaling is the standard
deviation of an entry of $W$ in units of $\sqrt{p_n}$. It tends to one when
$p_n\to0$, and then the scaling is $\sqrt{2d}$.

\subsection{Proof of Theorem~\ref{thm:semicircle}}

\begin{proof}[Proof of Theorem~\ref{thm:semicircle}]
	Let $d_{\mathrm{BL}}$ be the bounded Lipschitz metric~\eqref{eq:dbl}. By
	Lemma~\ref{lem:eigenvaluessemicircle},
	$d_{\mathrm{BL}}(\mu_n^\lambda,\mu_{sc})\to0$ in probability, and we must show the
	same for $\mu_n$. Every subsequence has a further subsequence along which
	$d_{\mathrm{BL}}(\mu_n^\lambda,\mu_{sc})\to0$ almost surely. It therefore suffices
	to show, for a deterministic sequence of graphs, that $\mu_n^\lambda\Rightarrow\mu_{sc}$
	implies $\mu_n\Rightarrow\mu_{sc}$.

	Let $X_k \defeq \tau_n(\lambda_k(L^tL) - 1)$ and
	$Y_k \defeq 2\tau_n(\sigma_k(L) - 1) = 2\tau_n\big(\sqrt{1+X_k/\tau_n} - 1\big)$.
	For $|x|\leq1/2$, $\sqrt{1 + x} = 1 + x/2 - x^2/(8(1+\xi)^{3/2})$ with $|\xi|\leq1/2$,
	and $(1+\xi)^{-3/2}\leq2\sqrt2$. Fix $M>2$. For $n$ so large that
	$M/\tau_n\leq1/2$, every $k$ with $|X_k|\leq M$ satisfies
	\begin{equation}\label{eq:taylorexpansion}
		\left|Y_k - X_k\right| \leq \frac{X_k^2}{\tau_n} \leq \frac{M^2}{\tau_n}.
	\end{equation}
	Fix $\epsilon>0$. Since $(\mu_n^\lambda)$ converges weakly, it is tight, and we
	choose $M$ with $\mu_n^\lambda(\{|x|>M\})<\epsilon$ for all large $n$. Let $f$ be
	in the test class of~\eqref{eq:dbl}, so $\Vert f\Vert_{\mathrm{Lip}}\leq1$ and
	$\Vert f\Vert_\infty\leq1$. Splitting the indices according to $|X_k|\leq M$,
	\begin{equation}
		\left|\mu_n^\lambda(f) - \mu_n(f)\right| \leq \frac{1}{n}\sum_{k=1}^n \left|f(X_k) - f(Y_k)\right| \leq \frac{M^2}{\tau_n} + 2\epsilon \leq 3\epsilon
	\end{equation}
	for all large $n$, uniformly in $f$. Hence
	$d_{\mathrm{BL}}(\mu_n^\lambda,\mu_n)\leq3\epsilon$ eventually, and
	$d_{\mathrm{BL}}(\mu_n,\mu_{sc})\to0$.
\end{proof}

\section{Discussion and outlook}\label{sec:discussion01}
\paragraph{Kemeny's constant.} By
Corollary~\ref{cor:averagedhittingtimes}, $H^i$ does not depend on $i$. It is
Kemeny's constant of the chain. Constancy has a generating-function proof with
relations for higher moments~\cite{fredes2025combinatorial}, the constant extends
to Markov processes~\cite{fitzsimmons2025kemeny}, and its spectral bounds and
computation are studied
in~\cite{abiad2025spectralkemeny,abiad2026kemenycompression,xia2024fast,xia2024means}.
Theorem~\ref{thm:hittingdirectedER} gives a law of large numbers for Kemeny's
constant of the simple random walk on $\overrightarrow{G}(n,p_n)$,
and~\eqref{eq:Hjfromp2p} the counterpart for the average starting hitting times.
For the undirected Erd\H{o}s--R\'enyi graph both satisfy central limit
theorems~\cite{lowe2025central,lowe_central_2023}.

\paragraph{The additive constant.}
Part~\ref{it:Kemenysecond} of Theorem~\ref{thm:hittingdirectedER} gives
$H^i = n-2+o(1)$ for $d\geq n^{\delta}$. Up to $o(1)$ this is Kemeny's constant
of the complete graph, $(n-1)^2/n = n-2+1/n$. In~\eqref{eq:kemenyintro}, the
first two terms of the trace expansion coincide with those of the complete graph,
and the higher terms vanish in probability.

\paragraph{What the estimates require of a model.} Formula~\eqref{eq:Hittingsingular}
holds for every chain. Its use on a random graph needs degree concentration and
control of the bulk of the singular values of the generator $L$. Bounds on the
non-Perron eigenvalues of $P$, such as those of~\cite[Theorem 2.12]{hazra2026directed}
for directed inhomogeneous random graphs, do not by themselves give the latter,
since $P$ is not normal. Nor do the singular values of $A$ determine those of $L$:
by Remark~\ref{rem:hermitianpart}, near the identity the latter follow the
Hermitian part of $P$.

\section*{Acknowledgements}
This work was partly funded by the Deutsche Forschungsgemeinschaft (DFG, German
Research Foundation), Project-ID 443891315, within SPP 2265. We gratefully
acknowledge the high-performance computing resources provided by the
Julius-Maximilians-Universit\"at W\"urzburg (JMU), funded by the DFG as a major
research instrumentation grant (\emph{Forschungsgro\ss ger\"at nach Art.~91b GG})
under INST 93/1145-1 FUGG.

\paragraph{Use of AI tools.} The main ideas, results, and original draft of this
paper were developed before generative AI tools became widely available. During
the final stage of manuscript preparation, the authors used large language
models to assist with language editing, computational checks, literature
searches, and internal critical review. The authors independently verified the
mathematical arguments, computations, references, and conclusions and take full
responsibility for the accuracy and content of the manuscript.

\appendix
\section{Standard tools}\label{app:tools}

\subsection{Matrix inequalities}\label{sec:svdprelim}

The following inequalities are standard; see~\cite[Chapter~3]{Horn_Johnson_1991}
and~\cite{bhatia2007perturbation}. For $A, B \in \mathbb{C}^{m \times n}$ and
$i,j\in\N$ with $i+j-1\leq\min\{m,n\}$, Weyl's inequalities state
\begin{equation}\label{eq:sumsingularvalues}
	\sigma_{i+j-1}(A + B) \leq \sigma_i(A) + \sigma_j(B),
\end{equation}
and, for $A,B\in\C^{n\times n}$~\cite[Theorem~3.3.16]{Horn_Johnson_1991},
$\sigma_{i+j-1}(AB)\leq\sigma_i(A)\sigma_j(B)$. Taking $i=1$, and applying the
same bound to $B = A^{-1}(AB)$, we obtain for invertible $A$ and every $k\in[n]$
\begin{equation}\label{eq:productsingularvalues}
	\sigma_n(A)\,\sigma_k(B) \leq \sigma_k(AB) \leq \sigma_1(A)\,\sigma_k(B).
\end{equation}
If $A,B\in\mathbb{C}^{n\times n}$ are Hermitian with eigenvalues
$\lambda_1(A)\leq\cdots\leq \lambda_n(A)$ and $\lambda_1(B)\leq\cdots\leq
\lambda_n(B)$, the Hoffman--Wielandt inequality~\cite{hoffman1953variation} gives
\begin{equation}\label{lidskii--wielandt:proposition}
	\sum_{i=1}^n \big(\lambda_i(A)-\lambda_i(B)\big)^2 \leq \Vert A-B\Vert_F^2.
\end{equation}

\subsection{Random matrices}\label{sec:rmtprelim}

For a Hermitian $n\times n$ matrix $A$, the empirical spectral distribution (ESD) is
$\mu_A \defeq \frac{1}{n}\sum_{i=1}^n \delta_{\lambda_i(A)}$. A sequence
$(W_n)_{n \in \mathbb{N}}$ of random Hermitian $n\times n$ matrices whose
upper-triangular entries are independent for each $n$ is a \emph{Wigner ensemble}.
The entries may depend on $n$. We use the triangular-array form of the semicircle
law~\cite[Theorem 2.9]{bai2010spectral}, which replaces moment conditions by a
Lindeberg condition.

\begin{theorem}[Semicircle law for Wigner ensembles]
	\label{thm:semicircleherensemble}
	Let $\left(W_n\right)_{n\in\N}$ be a Wigner ensemble whose upper-triangular
	entries $w_{ij}$ have mean zero and variance $\frac{1}{n}$. If for all
	$\epsilon > 0$
	\begin{equation}
		\lim_{n \to \infty} \frac{1}{n} \sum_{i,j=1}^n
		\mathbb{E}\left[ |w_{ij}|^2 \mathbbm{1}_{\{|w_{ij}| > \epsilon\}} \right]
		= 0,
	\end{equation}
	then $\mu_{W_n}$ converges weakly to the semicircle law, almost surely.
\end{theorem}

The bounded Lipschitz metric
\begin{align}
	\label{eq:dbl}
	d_{\mathrm{BL}}(\mu, \nu) \defeq \sup_{\substack{\Vert f\Vert_{\mathrm{Lip}}\leq 1\\ \Vert f\Vert_\infty\leq 1}} \left|\int f\,\mathrm{d}\mu-\int f\,\mathrm{d}\nu\right|
\end{align}
metrises weak convergence of probability measures on $\R$.

\begin{lemma}[ESD stability]
	\label{esd-stability:lemma}
	For Hermitian $n\times n$ matrices $A$ and $B$,
	$d_{\mathrm{BL}}(\mu_A,\mu_B) \leq n^{-1/2}\,\Vert A-B\Vert_F$.
\end{lemma}

\begin{proof}
	Order the eigenvalues non-decreasingly. For $f$ in the class of~\eqref{eq:dbl},
	by the Cauchy--Schwarz inequality and~\eqref{lidskii--wielandt:proposition},
	\begin{equation}
		\left|\int f\,\mathrm{d}\mu_A - \int f\,\mathrm{d}\mu_B \right| \leq \frac{1}{n}\sum_{i=1}^n \big|\lambda_i(A)-\lambda_i(B)\big|
		\leq \Big(\frac{1}{n}\sum_{i=1}^n \big(\lambda_i(A)-\lambda_i(B)\big)^2\Big)^{1/2}
		\leq \frac{\Vert A-B\Vert_F}{\sqrt{n}}. \qedhere
	\end{equation}
\end{proof}

\begin{corollary}[ESD convergence criterion]
	\label{cor:esd-stability}
	Let $B_n$ and $C_n$ be random Hermitian $n\times n$ matrices with
	$\frac1n\Vert B_n-C_n\Vert_F^2\to 0$ in probability. If $\mu_{B_n}$ converges
	weakly in probability to a probability measure $\mu$, then so does $\mu_{C_n}$.
\end{corollary}

\begin{proof}
	By Lemma~\ref{esd-stability:lemma}, $d_{\mathrm{BL}}(\mu_{B_n},\mu_{C_n})\to0$ in
	probability. Conclude by the triangle inequality for $d_{\mathrm{BL}}$.
\end{proof}

For Hermitian $n\times n$ matrices $A$ and $B$, the rank
inequality~\cite[Theorem~A.43]{bai2010spectral} bounds the Kolmogorov distance of
the distribution functions:
\begin{equation}\label{eq:rankinequality}
	\sup_{x\in\R}\left|\mu_A((-\infty,x]) - \mu_B((-\infty,x])\right| \leq \frac{\operatorname{rank}(A-B)}{n}.
\end{equation}

\begin{corollary}[ESD convergence under a finite-rank perturbation]
	\label{cor:esd-rank}
	Let $B_n$ and $C_n$ be random Hermitian $n\times n$ matrices with
	$C_n - B_n = K_n + E_n$, where $\operatorname{rank}(K_n)\leq r$ for a fixed $r$ and
	$\frac1n\Vert E_n\Vert_F^2\to 0$ in probability. If $\mu_{B_n}$ converges
	weakly in probability to a probability measure $\mu$ with a continuous
	distribution function, then so does $\mu_{C_n}$.
\end{corollary}

\begin{proof}
	Every subsequence has a further subsequence along which $\mu_{B_n}\Rightarrow\mu$
	almost surely. Since the distribution function of $\mu$ is continuous, weak
	convergence to $\mu$ is uniform convergence of the distribution functions.
	By~\eqref{eq:rankinequality} the latter passes from $B_n$ to $B_n+K_n$ with an
	error $r/n$. So $\mu_{B_n+K_n}\Rightarrow\mu$ in probability, and
	Corollary~\ref{cor:esd-stability} applied to $B_n+K_n$ and $C_n$ concludes.
\end{proof}

\begin{theorem}[Matrix Bernstein inequality~{\cite[Theorem 1.6]{tropp2012user}}]
	\label{thm:matrixbernstein}
	Let $X_1, \dots, X_N$ be independent random real $d_1 \times d_2$ matrices with
	$\E X_k = 0$ and $\Vert X_k\Vert \leq R$ almost surely. Let
	$S = \sum_{k=1}^N X_k$ and
	\begin{equation}
		\sigma^2 = \max \left\{
		\Big\Vert \sum_{k=1}^N \mathbb{E}[X_k X_k^t] \Big\Vert,\;
		\Big\Vert \sum_{k=1}^N \mathbb{E}[X_k^t X_k] \Big\Vert
		\right\}.
	\end{equation}
	Then, for all $t \ge 0$,
	\begin{equation}
		\P\left( \Vert S \Vert \ge t \right)\leq(d_1 + d_2) \exp\left( -\frac{ t^2 / 2 }{ \sigma^2 + Rt/3 }\right).
	\end{equation}
\end{theorem}

Finally, if $X$ is a sum of independent Bernoulli variables with mean $m$, then
for $0<t<1$~\cite[Corollary 4.6]{mitzenmacher2017probability}
\begin{equation}\label{eq:chernoff}
	\P\left(|X-m|\geq tm\right)\leq 2\exp\left(-\frac{t^2m}{3}\right).
\end{equation}

\section{Deferred proofs}\label{app:proofs}

\subsection{Concentration lemmas}\label{sec:concentrationproofs}

\begin{proof}[Proof of Lemma~\ref{lem:degreeconcentrationinout}]
	By the union bound it suffices to bound the probability that a fixed vertex,
	say $1$, violates the condition, and to multiply by $n$. The out-degree
	$d_1^+$ is a sum of $n-1$ independent Bernoulli$(p)$ variables, so its mean is
	$m \defeq (n-1)p = d(1 - 1/n)$. If $\left|d_1^+/d - 1\right| > \nu$, then
	\begin{equation}
		\left|d_1^+ - m\right| \geq \left|d_1^+ - d\right| - \left|d - m\right| > \nu d - \frac{d}{n} = t\,m, \qquad t \defeq \frac{\nu - 1/n}{1 - 1/n}\in(0,1).
	\end{equation}
	By~\eqref{eq:chernoff} the probability of this event is at most
	$2\exp(-t^2m/3)$. Since $n\geq10/\nu$,
	\begin{equation}
		t^2 m = \frac{d\left(\nu - 1/n\right)^2}{1 - 1/n} \geq d\left(\nu - \frac{\nu}{10}\right)^2 = \frac{81}{100}\,\nu^2 d \geq \frac34\,\nu^2d,
	\end{equation}
	which gives~\eqref{eq:degreeconcentrationout}. The in-degrees are treated in
	the same way.
\end{proof}

\begin{proof}[Proof of Lemma~\ref{lem:Wignernormbound}]
	Let $X \defeq W/\sqrt{d}$. Since $A_{ii}=0$,
	\begin{equation}\label{eq:Xsplit}
		X = \tilde{X} - \frac{p}{\sqrt{d}}I, \qquad \tilde{X} \defeq \sum_{i \neq j} Z^{(ij)},
		\qquad Z^{(ij)} \defeq \frac{A_{ij} - p}{\sqrt{d}} E^{(ij)},
	\end{equation}
	where $E^{(ij)}$ is the matrix unit at $(i,j)$. The $Z^{(ij)}$ are independent,
	centred, and $\Vert Z^{(ij)}\Vert \leq d^{-1/2} \eqdef R$. Moreover
	\begin{equation}
		\overrightarrow{\E}_{n,p}\big[Z^{(ij)}(Z^{(ij)})^t\big] = \frac{p(1-p)}{d}E^{(ii)} = \frac{1-p}{n}E^{(ii)},
	\end{equation}
	so $\sum_{i \neq j}\overrightarrow{\E}_{n,p}[Z^{(ij)}(Z^{(ij)})^t] = \frac{n-1}{n}(1-p) I$,
	and the same holds for $(Z^{(ij)})^tZ^{(ij)}$. Hence $\sigma^2\leq1$, and
	Theorem~\ref{thm:matrixbernstein} gives for $t > 0$
	\begin{equation}\label{eq:matrixBernstein}
		\overrightarrow{\P}_{n,p}\left( \Vert \tilde{X}\Vert \geq t \right) \leq 2n \exp\left( -\frac{t^2/2}{1 + \frac{t}{3\sqrt{d}}} \right).
	\end{equation}
	By~\eqref{eq:Xsplit}, $\Vert X\Vert \leq \Vert\tilde{X}\Vert + \sqrt{d}/n$, and
	$\sqrt{d}/n \leq \nu\sqrt{d}/2$ for $n \geq 2/\nu$. So
	$\{\Vert W\Vert > \nu d\}\subseteq\{\Vert\tilde{X}\Vert \geq \nu\sqrt{d}/2\}$, and
	$t = \nu\sqrt{d}/2$ in~\eqref{eq:matrixBernstein} gives~\eqref{eq:Wignernorm}.
\end{proof}

\subsection{Strong connectivity on the good events}\label{sec:connectivityproof}

\begin{lemma}[Strong connectivity]\label{lem:stronglyconnected}
	Let $\nu_1,\nu_2>0$ with $\nu_1+\nu_2<1/2$. On $\cF_{\nu_1}\cap\cH_{\nu_2}$ the
	digraph is strongly connected.
\end{lemma}

\begin{proof}
	Let $S\subseteq V$ be nonempty with no arc from $S$ to $V\setminus S$, and let
	$k\defeq|S|$. Every arc leaving a vertex of $S$ ends in $S$, so by $\cF^+_{\nu_1}$
	\begin{equation}
		k(1-\nu_1)d \leq \sum_{i\in S}d_i^+ = \mathbf{1}_S^tA\mathbf{1}_S
		= pk^2 + \mathbf{1}_S^tW\mathbf{1}_S \leq pk^2 + \nu_2dk,
	\end{equation}
	using $\mathbf{1}_S^tJ\mathbf{1}_S = k^2$ and $\cH_{\nu_2}$. Since $pk = dk/n$, this
	gives $k\geq(1-\nu_1-\nu_2)n>n/2$. The same argument with in-degrees, using
	$\cF^-_{\nu_1}$ and $\Vert W^t\Vert = \Vert W\Vert$, shows that every nonempty
	$T\subseteq V$ with no arc from $V\setminus T$ to $T$ has $|T|>n/2$.

	If the digraph is not strongly connected, its strongly connected components form
	an acyclic digraph with at least two vertices. A sink component $S$ has no
	outgoing arcs and a source component $T$ has no incoming arcs, and $S\cap
	T=\emptyset$. Then $|S|+|T|>n$, a contradiction.
\end{proof}

\subsection{Proof of Lemma~\ref{lem:Bsemicircle}}

The off-diagonal entries of $W = A-pJ$ are independent, so the entries $b_{ij}$ of $B_n$ with
$i<j$ are independent. Fix $c<1$ with $p\leq c$ for $n$ large. For $i\neq j$,
$b_{ij} = -(W_{ij}+W_{ji})/\sqrt{2d(1-p)}$ has mean zero, variance
$2p(1-p)/(2d(1-p)) = 1/n$, and $|b_{ij}|\leq\sqrt{2/(d(1-c))}$. The diagonal is
deterministic: $b_{ii} = 2p/\sqrt{2d(1-p)} \eqdef c_n$, and $c_n\leq\sqrt{2/(n(1-c))}$
since $d\leq n$. The theorem we cite asks for centred diagonal entries of variance $1/n$. To match
this normalisation only, let $D_n$ be a diagonal matrix of independent signs,
scaled by $n^{-1/2}$ and independent of $A$, and set
\begin{equation}
	\tilde B_n \defeq B_n - c_nI + D_n.
\end{equation}
The upper-triangular entries of $\tilde B_n$ are independent, centred, with
variance $1/n$, and in absolute value at most the larger of $\sqrt{2/(d(1-c))}$
and $n^{-1/2}$, which tends to $0$. So the
Lindeberg sum of Theorem~\ref{thm:semicircleherensemble} vanishes for $n$ large,
and $\mu_{\tilde B_n}$ converges weakly to the semicircle law almost surely. By
Lemma~\ref{esd-stability:lemma},
\begin{equation*}
	d_{\mathrm{BL}}\bigl(\mu_{\tilde B_n},\mu_{\tilde B_n-D_n}\bigr)\leq n^{-1/2}\Vert D_n\Vert_F = n^{-1/2},
\end{equation*}
so the same holds for
$\tilde B_n - D_n = B_n - c_nI$. Finally, $\mu_{B_n}$ is the image of
$\mu_{B_n-c_nI}$ under $x\mapsto x+c_n$, and $c_n\to0$. \qed

\section{Alternative proof of part~\ref{it:Hitarget} via singular values}\label{sec:svdroute}

This appendix is not needed for the main results. We give a second proof of part~\ref{it:Hitarget} of
Theorem~\ref{thm:hittingdirectedER}, through the formula of
Theorem~\ref{thm:Hittingsingular}. It needs the location of the spectral edges of
$L$ (Lemma~\ref{lem:singularbulkconvergence}), a bound on the sum of the inner
products $v_k^tu_k$ (Lemma~\ref{lem:principalangles}), and a bound on an error term
(Lemma~\ref{lem:errorbounds}). Lemma~\ref{lem:statalignment} is shared with the
main proof.

\begin{lemma}[Concentration of the singular values of $L$]
	\label{lem:singularbulkconvergence}
	For every $\epsilon > 0$ there exist $\nu_1\in(0,1)$ and $\nu_2>0$ such that, on
	$\cF_{\nu_1}\cap\cH_{\nu_2}$,
	\begin{equation}
		1 - \epsilon < \sigma_{n-1}(L) \leq \sigma_1(L) < 1 + \epsilon.
	\end{equation}
\end{lemma}

\begin{proof}
	Since $A = pJ + W$, we have $L = D_+^{-1}(D_+ - pJ - W)$. By~\eqref{eq:productsingularvalues},
	\begin{equation}
		\sigma_{n-1}(L) \geq \sigma_n(D_{+}^{-1})\,\sigma_{n-1}\left(D_+ - pJ - W\right),
		\qquad
		\sigma_{1}(L) \leq \sigma_1(D_{+}^{-1})\,\sigma_{1}\left(D_+ - pJ - W\right).
	\end{equation}
	By~\eqref{eq:sumsingularvalues}, applied twice,
	\begin{equation}
		\sigma_{n-1}\left(D_+ - pJ - W\right) \geq \sigma_{n-1}\left(D_+ - pJ\right) - \Vert W\Vert \geq \sigma_n\left(D_+ \right) - \Vert W\Vert.
	\end{equation}
	The second inequality is Weyl's inequality with $D_+=(D_+-pJ)+pJ$,
	$i=n-1$, $j=2$, and $\sigma_2(pJ)=0$. The matrix $pJ-dI$ has eigenvalues
	$0$ and $-d$, so $\Vert pJ-dI\Vert = d$, and
	\begin{equation}
		\sigma_1\left(D_+ - pJ - W\right) \leq \Vert D_+ - dI\Vert + \Vert W\Vert + \Vert pJ - dI\Vert.
	\end{equation}
	On $\cF_{\nu_1}\cap\cH_{\nu_2}$ every out-degree lies in $[d(1-\nu_1),d(1+\nu_1)]$
	and $\Vert W\Vert\leq\nu_2d$. Hence
	\begin{equation}
		\frac{1 - \nu_1 - \nu_2}{1 + \nu_1} \leq \sigma_{n-1}(L) \leq \sigma_1(L) \leq \frac{1 + \nu_1 + \nu_2}{1-\nu_1}.
	\end{equation}
	Both bounds tend to $1$ as $\nu_1,\nu_2\to0$.
\end{proof}

Near-normality of $L$ does not control the inner products $v_k^tu_k$. If
$L = \sum_{k=1}^n\lambda_kw_kw_k^*$ is normal with $\lambda_n = 0$ simple, its
singular triplets for $k<n$ are $(|\lambda_k|, w_k, \frac{\lambda_k}{|\lambda_k|}w_k)$, so
\begin{equation}\label{eq:normalphasesum}
	\sum_{k=1}^{n-1}v_k^*u_k = \sum_{k=1}^{n-1}\frac{\lambda_k}{|\lambda_k|}.
\end{equation}
For a reversible chain, after the similarity transform~\eqref{eq:Lhat}, all
phases equal $1$ and the sum is $n-1$. For the oriented cycle of
Example~\ref{ex:MCtotheright}, $L$ is circulant, hence normal, but
$\lambda_k = 1 - e^{2\pi\mathrm{i}k/n}$ has phase $\sin(\pi k/n) - \mathrm{i}\cos(\pi k/n)$, and
\begin{equation}\label{eq:cyclephasesum}
	\sum_{k=1}^{n-1}\frac{\lambda_k}{|\lambda_k|} = \sum_{k=1}^{n-1}\sin\left(\frac{\pi k}{n}\right) = \cot\left(\frac{\pi}{2n}\right) = \left(\frac{2}{\pi} + o(1)\right)n.
\end{equation}
So no bound on $LL^t - L^tL$ can give $\sum_kv_k^tu_k = n(1+o(1))$. The control comes instead
from the Frobenius distance $\Vert P-\frac1nJ\Vert_F$, through the next lemma.

\begin{lemma}
	\label{lem:principalangles}
	Let $P$ be irreducible and let $Q \defeq P - \frac{1}{n}J$. Then
	\begin{equation}\label{eq:principalanglesexact}
		0 \leq (n-1) - \sum_{k=1}^{n-1} v_k^tu_k \leq 2\sum_{k=1}^{n-1}\left\Vert Pv_k\right\Vert^2 = 2\left\Vert Q\right\Vert_F^2.
	\end{equation}
	In particular, on $\cF^+_{\nu_1}$,
	\begin{equation}\label{eq:principalanglesbulk}
		0 \leq (n-1) - \sum_{k=1}^{n-1} v_k^tu_k \leq \frac{2n}{\left(1-\nu_1\right)d} + 2.
	\end{equation}
\end{lemma}

\begin{proof}
	Let $k<n$. Then $\sigma_k>0$, since $P$ is irreducible and $\sigma_n(L) = 0$ is
	simple. Write $w_k \defeq Lv_k = v_k - Pv_k$, so that $u_k = w_k/\Vert w_k\Vert$.
	The sign of $v_k$ is free, but it fixes the sign of $u_k$, and $v_k^tu_k$ is
	unchanged when both flip.
	Both $u_k$ and $v_k$ are unit vectors, so $1 - v_k^tu_k = \frac12\Vert u_k-v_k\Vert^2\geq 0$.
	For a unit vector $v$ and $w\neq0$,
	\begin{equation}
		\left\Vert \frac{w}{\Vert w\Vert} - v\right\Vert \leq \left\Vert\frac{w}{\Vert w\Vert} - w\right\Vert + \Vert w - v\Vert = \bigl|1 - \Vert w\Vert\bigr| + \Vert w-v\Vert \leq 2\Vert w - v\Vert.
	\end{equation}
	With $w = w_k$ and $v = v_k$ this gives $\Vert u_k - v_k\Vert\leq 2\Vert Pv_k\Vert$,
	hence $0\leq 1 - v_k^tu_k\leq 2\Vert Pv_k\Vert^2$. Summing over $k<n$ gives the
	inequalities in~\eqref{eq:principalanglesexact}.

	Since $L\mathbf{1}_n = 0$, we take $v_n = n^{-1/2}\mathbf{1}_n$. Then
	$(v_k)_{k\in[n-1]}$ is an orthonormal basis of $\mathbf{1}_n^\perp$, and
	$Pv_k = Qv_k$ for $k<n$, since $Jv_k = 0$. Also $Q\mathbf{1}_n = 0$, so
	$\sum_{k<n}\Vert Pv_k\Vert^2 = \sum_{k\leq n}\Vert Qv_k\Vert^2 = \Vert Q\Vert_F^2$.

	On $\cF^+_{\nu_1}$ every $d_i^+\geq(1-\nu_1)d$. Each entry of row $i$ of $Q$ with
	$A_{ij}=1$ has square at most $(d_i^+)^{-2}$, since $1/n\leq1/d_i^+$, and each
	other entry has square $n^{-2}$. Summing,
	\begin{equation}\label{eq:QFrobenius}
		\left\Vert Q\right\Vert_F^2 \leq \sum_{i\in[n]}\frac{1}{d_i^+} + 1 \leq \frac{n}{\left(1-\nu_1\right)d} + 1,
	\end{equation}
	and~\eqref{eq:principalanglesbulk} follows.
\end{proof}

The right-hand side of~\eqref{eq:principalanglesbulk} is $O(n/d) = o(n)$. It uses
only the out-degrees. For the oriented cycle $\Vert Q\Vert_F^2 = n-1$,
and~\eqref{eq:principalanglesexact} gives no information, in accordance
with~\eqref{eq:cyclephasesum}.

\begin{lemma}[Error term]\label{lem:errorbounds}
	Let $P$ be irreducible with $\sigma_{n-1}(L) > 0$, and let
	$x = \pi - \frac{1}{n}\mathbf{1}_n$. Then, for every $k \in [n-1]$,
	\begin{align}
		\left|v_k^t\pi\right| \leq \Vert x\Vert, \qquad
		\left|u_k^t\mathbf{1}_n\right| \leq n\Vert x\Vert, \qquad
		\left|\frac{1}{n}\sum_{k=1}^{n-1} \frac{1}{\sigma_k}v_k^t\pi\, u_k^t\mathbf{1}_n\right| \leq \frac{n\Vert x\Vert^2}{\sigma_{n-1}}. \label{eq:hittingsingularerrorup}
	\end{align}
\end{lemma}

\begin{proof}
	We have $v_n \propto \mathbf{1}_n$ and $u_n\propto\pi$, so $v_k^t\mathbf{1}_n = 0$
	and $u_k^t\pi = 0$ for $k<n$. Hence $v_k^t\pi = v_k^tx$ and
	$u_k^t\mathbf{1}_n = -n\,u_k^tx$, and the Cauchy--Schwarz inequality gives the
	first two bounds. The third follows by summing.
\end{proof}

\begin{proof}[Second proof of part~\ref{it:Hitarget} of Theorem~\ref{thm:hittingdirectedER}]
	Fix $\epsilon' > 0$. By Lemmas~\ref{lem:singularbulkconvergence}
	and~\ref{lem:statalignment} there are $\nu_1, \nu_2 > 0$ such that, on
	$\cF_{\nu_1}\cap\cH_{\nu_2}$ and for $n$ large,
	\begin{equation}
		1-\epsilon' < \sigma_{n-1}(L) \leq \sigma_1(L) < 1+\epsilon' \qquad\text{and}\qquad n\Vert x\Vert^2 < \epsilon'.
	\end{equation}
	By~\eqref{eq:meantargethittingtime},
	$H^i = \sum_{k=1}^{n-1} \sigma_k^{-1}\left(v_k^tu_k - v_k^t\pi\, u_k^t\mathbf{1}_n\right)$.
	By~\eqref{eq:hittingsingularerrorup}, the second part contributes at most
	$\epsilon'/(1-\epsilon')$ to $H^i/n$. For the first part let
	$a_k \defeq v_k^tu_k\leq1$ and $D \defeq \sum_{k=1}^{n-1}(1-a_k)\geq0$. Since
	$\sum_{k:\,a_k\geq0} a_k \leq \sum_k a_k + \sum_{k:\,a_k<0}(1-a_k)\leq \sum_ka_k+D = n-1$,
	\begin{equation}\label{eq:hittingsingularmain}
		\frac{n-1}{n\,\sigma_1(L)} - \frac{D}{n\,\sigma_{n-1}(L)}
		\leq \frac{1}{n}\sum_{k=1}^{n-1} \frac{a_k}{\sigma_k}
		\leq \frac{1}{n\,\sigma_{n-1}(L)}\sum_{k:\,a_k\geq0} a_k \leq \frac{1}{\sigma_{n-1}(L)},
	\end{equation}
	the lower bound from $\sum_k a_k/\sigma_k = \sum_k1/\sigma_k - \sum_k(1-a_k)/\sigma_k$.
	By Lemma~\ref{lem:principalangles},
	\begin{equation}\label{eq:Dsmall}
		\frac{D}{n} \leq \frac{2}{\left(1-\nu_1\right)d} + \frac{2}{n} < \epsilon'
	\end{equation}
	for $n$ large, since $d\to\infty$. So $|H^i/n - 1|<\epsilon$ on
	$\cF_{\nu_1}\cap\cH_{\nu_2}$ for $n$ large, once $\epsilon'$ is small enough.
	Corollary~\ref{cor:goodprobability} then gives a bound of the same form
	as~\eqref{eq:meanhittingdirectedER}, with the constants adjusted.
\end{proof}

\section{The oriented cycle in closed form}\label{sec:circulantexample}

We compute the singular triplets of the generator of Example~\ref{ex:MCtotheright}
and evaluate Corollary~\ref{cor:averagedhittingtimes} on them. Let $n$ be even; the
odd case is similar. We index the triplets by the Fourier frequency $r\in[n]$, not by $k$, so
that the triplet $r$ carries the singular value $2\sin(\pi(r-1)/n)$. This differs
from the non-increasing order of Section~\ref{sec:notation} by a permutation, which
does not affect the sums in Theorem~\ref{thm:Hittingsingular} and
Corollary~\ref{cor:averagedhittingtimes}.

Let $c_r = n^{-1/2}$ for $r = n/2+1$ (the Nyquist mode) and $c_r = (2/n)^{1/2}$
otherwise:
\begin{equation}\label{eq:circulantnormalisation}
	c_r^2 = \frac{2}{n} - \frac1n\mathbbm{1}_{\{r=n/2+1\}}.
\end{equation}
Let $v_1 = u_1 = n^{-1/2}\mathbf{1}_n$, and for $j\in[n]$ and $\theta_r\defeq\pi(r-1)/n$,
\begin{equation}
	(v_r)_j = c_r\cos\left(2\theta_r(j-1)\right),\quad (u_r)_j = c_r\sin\left(\theta_r(2j-1)\right)\qquad (2\leq r\leq n/2+1),
\end{equation}
\begin{equation}
	(v_r)_j = c_r\sin\left(2\theta_r(j-1)\right),\quad (u_r)_j = -c_r\cos\left(\theta_r(2j-1)\right)\qquad (n/2+2\leq r\leq n).
\end{equation}
Both families are orthonormal by the product-to-sum formulas. Since
$((I-P)v)_j = v_j - v_{j+1}$, the same formulas give $(I-P)v_r = 2\sin(\theta_r)\,u_r$
for every $r$. Thus $(2\sin\theta_r, v_r, u_r)$ are singular triplets of $L$.

Here $\pi = n^{-1}\mathbf{1}_n$, so $v_r^t\pi = 0$ and $u_r^t\mathbf{1}_n = 0$ for
$r\geq2$, and $H^i = \sum_{r=2}^n v_r^tu_r/(2\sin\theta_r)$. For $2\leq r\leq n/2+1$,
\begin{equation}
	v_r^tu_r = \frac{c_r^2}{2} \sum_{j=1}^n \left[\sin\left(\theta_r(4j-3)\right) + \sin\theta_r\right].
\end{equation}
The first sum vanishes for $r\leq n/2$ and equals $n$ for $r = n/2+1$. In both cases
$v_r^tu_r = \sin\theta_r$, by~\eqref{eq:circulantnormalisation}. The sine branch is
identical. Hence
\begin{equation}\label{eq:meantargethittingtimecirculant}
	H^i = \sum_{r=2}^n \frac{\sin\theta_r}{2\sin\theta_r} = \frac{n-1}{2}.
\end{equation}
For the average starting hitting time at $j=1$, $(v_r)_1(u_r)_1 = c_r^2\sin\theta_r$
for $r\leq n/2+1$ and $0$ otherwise, so
\begin{equation}
	H_1 = n\sum_{r=2}^{n/2+1} \frac{c_r^2 \sin\theta_r}{2\sin\theta_r} = \frac{n}{2}\left(\frac{2}{n}\left(\frac n2-1\right) + \frac1n\right) = \frac{n-1}{2}.
\end{equation}
By vertex-transitivity, $H_j = (n-1)/2$ for every $j$.

\printbibliography

\end{document}